\documentclass{article}
\usepackage[utf8]{inputenc}

\usepackage{amsmath, amssymb, amsthm, mathtools}
\usepackage{pgfplots}
\pgfplotsset{compat=1.14}
\usepackage{tikz}
\usetikzlibrary{shapes}
\usepackage{tabularx}
\usepackage{multirow}
\usepackage{float}
\usepackage[usestackEOL]{stackengine} 

\newtheorem{Thm}{Theorem}[section]

\newtheorem{Prb}[Thm]{Problem}
\theoremstyle{definition}
\newtheorem{Def}[Thm]{Definition}
\newtheorem{question}[Thm]{Question}

\newcommand{\N}{\mathbb{N}}
\newcommand{\Z}{\mathbb{Z}}
\newcommand{\R}{\mathbb{R}}

\newcommand{\T}{\mathbb{T}}
\newcommand{\mmin}[1]{\min \left\{ #1 \right\}}
\newcommand{\mmax}[1]{\max \left\{ #1 \right\}}
\newcommand{\convexhull}[1]{\operatorname{conv}\left(#1\right)}
\newcommand{\GL}[1]{\operatorname{GL_{#1}}(\Z)}
\newcommand{\of}[1]{\left(#1\right)}
\newcommand{\off}[1]{\left[#1\right]}
\newcommand{\offf}[1]{\left\{#1\right\}}
\newcommand{\st}[1]{\; \left| \; #1 \right.}
\newcommand{\abs}[1]{\lvert #1 \rvert}
\newcommand{\triang}[4]{\begin{tikzpicture}[scale = 0.25] \draw (0,0)--(#1,#2)--(#3,#4)--cycle; \end{tikzpicture}}
\newcommand{\standardsimplex}{\begin{tikzpicture}[scale = 0.22] \draw (0,0)--(1,0)--(0,1)--cycle; \end{tikzpicture}}
\newcommand{\standardtetrahedron}{
    \begin{tikzpicture}[scale = 0.22]
        \coordinate (O) at (0,0,0);
        \coordinate (A) at (1,0,0);
        \coordinate (B) at (0,1,0);
        \coordinate (C) at (0,0,1);
        \draw (O)--(A);
        \draw (O)--(B);
        \draw (O)--(C);
        \draw (A)--(B)--(C)--cycle;
    \end{tikzpicture}
}
\newcommand{\vertices}[9]{                                                     
\coordinate (O) at (0,0,0);                                                     
\coordinate (A) at (#1,#2,#3);                                                  
\coordinate (B) at (#4,#5,#6);                                                  
\coordinate (C) at (#7,#8,#9);                                                  
}                                                                                
\newcommand{\edges}[1]{                                                        
\filldraw[fill=#1, opacity=0.5] (O)--(A)--(B)--cycle;                           
\filldraw[fill=#1, opacity=0.5] (O)--(A)--(C)--cycle;                           
\filldraw[fill=#1, opacity=0.5] (O)--(B)--(C)--cycle;                           
\filldraw[fill=#1, opacity=0.5] (A)--(B)--(C)--cycle;                           
} 
\renewcommand{\vector}[2]{\left( \begin{smallmatrix}#1\\#2\end{smallmatrix} \right)}

\newcommand{\conv}[6]{\convexhull{\offf{\vector{#1}{#2},\vector{#3}{#4},\vector{#5}{#6}}}}

\newcolumntype{C}{>{$} c <{$}}
\newcolumntype{L}{>{$} l <{$}}
\newcolumntype{R}{>{$} r <{$}}
\newcolumntype{P}[1]{>{\centering\arraybackslash $} p{#1} <{$}}
\newcolumntype{Q}[1]{>{\raggedright\arraybackslash $} p{#1} <{$}}

\begin{document}

\title{Algorithmic Symplectic Packing}
\author{Greta Fischer\and Jean Gutt \and Michael J{\"u}nger}

\maketitle

\begin{abstract}
\noindent In this article we explore a symplectic packing problem where the targets and domains are $2n$-dimensional symplectic manifolds. We work in the context where the manifolds have first homology group equal to~$\Z^n$, and we require the embeddings to induce isomorphisms between first homology groups. 
In this case, Maley, Mastrangeli and Traynor~\cite{Maley2000} showed that the problem can be reduced to a combinatorial optimization problem, namely packing certain allowable simplices into a given standard simplex.
They designed a computer program and presented computational results. In particular, they determined the simplex packing widths in dimension four for up to $k=12$ simplices, along with lower bounds for higher values of $k$.
We present a modified algorithmic approach
that allows us to determine the 
$k$-simplex packing widths for up to $k = 13$ simplices in dimension four and up to $k = 8$ simplices in dimension six.
Moreover, our approach determines all simplex-multisets that allow for optimal packings.


\end{abstract}

\section*{Introduction and Motivation}

Symplectic geometry originated in a mathematical formulation of the classical mechanics of dynamical systems with finitely many degrees of freedom. The objects studied are smooth manifolds with an additional structure. Formally, a symplectic manifold ($M,\omega$) is a pair consisting of a smooth manifold $M$ and a symplectic form $\omega$ on $M$. A symplectic form $\omega$ on a manifold $M$ is a closed non-degenerate 2-form on $M$. A map $\varphi: (M_1, \omega_1) \rightarrow (M_2, \omega_2)$ between two symplectic manifolds is called symplectic if it satisfies $\varphi ^\star \omega_2 = \omega_1$. As Riemannian geometry is the study of transformations preserving the inner product, symplectic geometry is  the study  of  transformations  preserving  the symplectic form. An important theorem in symplectic geometry, which goes back to Darboux in 1882, is that locally all symplectic forms are the same. This result implies that there are no local invariants in symplectic geometry. This is in great contrast to Riemannian geometry where curvature is an invariant that can be determined locally. So the natural question that arises is: Can we find globally defined symplectic invariants? 

One of the first striking results in this direction, which lies at the root of symplectic geometry, is due to Gromov. He asked for the greatest radius of a ball that can be symplectically embedded into a given symplectic manifold. To illustrate the nontriviality of this question, he stated the famous Non-squeezing Theorem in 1985~\cite{Gromov1985}. Let 
\begin{align*}
    B^{2n}(r) &= \offf{(x, y)\in \R^{2n} \st{\abs{x}^2 + \abs{y}^2 < \frac{r}{\pi}}},\\
    Z^{2n}(s) &= \offf{(x, y)\in \R^{2n} \st{x_1 ^2 + y_1^2 < \frac{s}{\pi}}}
\end{align*}
denote the $2n$-dimensional open ball of radius $\sqrt{\frac{r}{\pi}}$ and the $2n$-dimensional open cylinder of radius $\sqrt{\frac{s}{\pi}}$, respectively. Gromov's Non-squeezing theorem states that one cannot symplectically embed $B^{2n}(r)$ into $Z^{2n}(s)$ unless the radius $r$ of the ball is less than or equal to the radius $s$ of the cylinder. The Non-squeezing Theorem indicates the rigidity of symplectic embeddings as compared to the flexibility of volume preserving diffeomorphisms. 

The invariant found by Gromov is called the ball packing width. Let the expression $\varphi: (M_1, \omega_1) \xhookrightarrow{s} (M_2, \omega_2)$ denote that the map $\varphi$ is a symplectic embedding. Then the ball packing width is formally defined as follows. 
\begin{Def}
    The ball packing width of a $2n$-dimensional symplectic manifold $(M, \omega)$ is 
    \begin{equation*}
        g(M, \omega) = \operatorname{sup} \offf{r \st{\exists \varphi: \of{B^{2n}(r), \omega_0} \xhookrightarrow{s} (M, \omega)}}.
    \end{equation*}
\end{Def} 
Instead of studying only one symplectic embedding of a ball of maximum radius, one can also study $k$ symplectic embeddings of a ball with maximum radius such that the embeddings have pairwise disjoint images. The corresponding invariant is called the $k$-ball packing width.
\begin{Def}
    The $k$-ball packing width of a $2n$-dimensional symplectic manifold $(M, \omega)$ is
    \begin{equation*}
        g_k(M, \omega) = \operatorname{sup} \offf{r\st{
        \begin{array}{l}
        \exists \varphi_1,\ldots,\varphi_k: \of{B^{2n}(r), \omega_0} \xhookrightarrow{s} (M, \omega) \text{ with } \\
        \varphi_i\of{B^{2n}(r)} \cap \varphi_j\of{B^{2n}(r)} = \emptyset \quad \forall 1\leq i < j \leq k    
        \end{array}
        }}.
    \end{equation*}
\end{Def} 
The ball packing width and the $k$-ball packing width are symplectic invariants which have been studied in~\cite{Biran1996, Maley2000, McDuff1994, Schlenk2005, Traynor1995, Wieck2009}. They give information that can be used to distinguish symplectic manifolds. 
Key properties of these and other symplectic invariants were first axiomatized in 1994 by Ekeland and Hofer who introduced the notion of symplectic capacity~\cite{Hofer2012}.
They define a symplectic capacity as a map $c$ from the set of all symplectic manifolds to the interval $\off{0, \infty}$ that satisfies the following three properties:
\begin{enumerate}
    \item Monotonicity: $(M_1, \omega_1) \xhookrightarrow{s} (M_2, \omega_2) \Rightarrow c(M, \omega_1) \leq c(M_2, \omega_2)$.
    \item Conformality: $\forall \alpha \in \R \setminus \{0\}: \; c(M, \alpha \omega) = |\alpha| c(M, \omega)$.
    \item Nontriviality: $c \of{B^{2n}(1), \omega_0} > 0$ and $c \of{Z^{2n}(1), \omega_0} < \infty$.
\end{enumerate}
The nontriviality condition guarantees that, in dimension greater than two, volume is not a capacity. The search for symplectic capacities and techniques to calculate them are major areas of research in symplectic geometry. Although these invariants are quite easy to define, they are extremely difficult to calculate. 
For a computational approach, we replace the ball by an open prism approximation.
Let
\begin{equation*}
    P^{2n}(r) = \T^n \times \standardsimplex^n(r)
\end{equation*}
denote the $2n$-dimensional open prism, which is the Cartesian product of the $n$-dimensional torus $\T^n = \R^n / \Z^n$ and the $n$-dimensional open standard simplex of side length $r$
\begin{equation*}
    \standardsimplex^n(r) = \Big\{x \in \R^{n} \;\Big\vert\; x_i > 0 \quad \forall i \in \offf{1,\ldots,n} \text{ and } \sum_{i=1}^{n} x_i < r \Big\}.    
\end{equation*}
We endow $P^{2n}(r)$ with the symplectic form $\omega_0$ coming from the standard symplectic form on $\R^{2n}$.
Mastrangeli~\cite{Mastrangeli1997} has shown that, for every $\varepsilon > 0$, there exist symplectic embeddings 
    \begin{equation*}
        \of{B^{2n}(r-\varepsilon),\omega_0} \xhookrightarrow{s} \of{P^{2n}(r),\omega_0} \xhookrightarrow{s} \of{B^{2n}(r),\omega_0}.
    \end{equation*}
Therefore, embeddings of the ball into a symplectic manifold give rise to embeddings of the prism and vice versa. For this reason, there is no quantitative difference between looking at symplectic packings of the ball or the prism. 
When we restrict to symplectic manifolds that have first homology equal to $\Z^n$, it is possible to add the condition that the symplectic embeddings induce isomorphisms on the level of first homology. We call these maps $1$-isomorphic. With this additional property we can define a slight modification of the $k$-ball packing width that allows for a computational approach.
\begin{Def}
    The $k$-simplex packing width of a $2n$-dimensional symplectic manifold $(M, \omega)$ is 
    \begin{equation*}
        s_k(M, \omega) = \operatorname{sup} \offf{r\st{
        \begin{array}{l}
        \exists \varphi_1,\ldots,\varphi_k: \of{P^{2n}(r), \omega_0} \xhookrightarrow[1-isomorphic]{s} (M, \omega) \text{ with } \\
        \varphi_i\of{P^{2n}(r)} \cap \varphi_j\of{P^{2n}(r)} = \emptyset \quad \forall 1\leq i < j \leq k    
        \end{array}
        }}\!.
    \end{equation*}
\end{Def} 

The $k$-simplex packing width has been introduced by Maley et~al.~\cite{Maley2000}.
It is a symplectic invariant that satisfies a set of axioms analogous to the capacity axioms.  Since we reduce the set of symplectic embeddings to those that are $1$-isomorphic, the $k$-simplex packing width gives a lower bound on the $k$-ball packing width: $s_k \leq g_k$. 

Maley et~al.\ calculated the $k$-simplex packing width of $P^4(1)$ by using a computer program. Their results are shown in Table~\ref{tab: packing-widths}. The $k$-ball packing width can be computed using pseudoholomorphic curves~\cite{McDuff1994}. For $k > 12$ the values of $s_k \of{P^4(1), \omega_0}$ are known lower bounds that are conjectured to be optimal. The approach of Maley et~al.\ builds up on work of Traynor who obtained maximal packings of $B^4(1)$ by constructing maximal packings of $P^4(1)$ by hand~\cite{Traynor1995}. Similar techniques were applied by Kruglikov~\cite{Kruglikov1996}. Other explicit symplectic packings were constructed by Schlenk~\cite{Schlenk2005} and Wieck~\cite{Wieck2009}.

\begin{table}[ht]
    \renewcommand{\arraystretch}{1.2}
    \begin{small}
    \begin{center}
    {\tabcolsep 1.3mm%
    \begin{tabular}{|P{0.75cm}|P{0.75cm}P{0.75cm}P{0.75cm}P{0.75cm}P{0.75cm}P{0.75cm}P{0.75cm}P{0.75cm}P{0.75cm}P{0.75cm}|}
            \hline
            k & 1 & 2 & 3 & 4 & 5 & 6 & 7 & 8 & 9 & 10\\ 
            \hline 
            g_k & 1 & \frac{1}{2} & \frac{1}{2} & \frac{1}{2} & \frac{2}{5} & \frac{2}{5} & \frac{3}{8} & \frac{6}{17} & \frac{1}{3} & \frac{1}{\sqrt{10}}\\ 
            s_k & 1 & \frac{1}{2} & \frac{1}{2} & \frac{1}{2} & \frac{2}{5} & \frac{6}{17} & \frac{1}{3} & \frac{1}{3} & \frac{1}{3} & \frac{3}{10} \\ 
            \hline
        \end{tabular}
        \\[5mm]
        \begin{tabular}{|P{0.75cm}|P{0.75cm}P{0.75cm}P{0.75cm}P{0.75cm}P{0.75cm}P{0.75cm}P{0.75cm}P{0.75cm}P{0.75cm}P{0.75cm}|}
        \hline
        k & 11 & 12 & 13 & 14 & 15 & 16 & 17 & 18 & 19 & 20 \\ %& 21 & 22 & 23 \\ 
        \hline 
        g_k & \frac{1}{\sqrt{11}} & \frac{1}{\sqrt{12}} & \frac{1}{\sqrt{13}} & \frac{1}{\sqrt{14}} & \frac{1}{\sqrt{15}} & \frac{1}{4} & \frac{1}{\sqrt{17}} & \frac{1}{\sqrt{18}} & \frac{1}{\sqrt{19}} & \frac{1}{\sqrt{20}} \\ %& \frac{1}{\sqrt{21}} & \frac{1}{\sqrt{22}} & \frac{1}{\sqrt{23}} \\
        s_k & \frac{2}{7} & \frac{15}{56} & \geq \frac{6}{23} & \geq \frac{20}{79} & \geq \frac{1}{4}  & \frac{1}{4} & \geq \frac{4}{17} & \geq \frac{3}{13} & \geq \frac{2}{9} & \geq \frac{21}{97} \\ %& \geq \frac{4}{19} & \geq \frac{7}{34} & \geq \frac{21}{104} \\
        \hline
    \end{tabular}
    }
    \end{center}
    \end{small}
    \caption{$g_k \of{P^4(1), \omega_0}$ versus $s_k \of{P^4(1), \omega_0}$~\cite{Maley2000}}
    \label{tab: packing-widths}
\end{table}

The implementation of the algorithmic approach that we are going to present confirms the results of~\cite{Maley2000} for $k$ up to 12 and confirms the optimality of the lower bound for the simplex packing width for $k=13$.

Moreover, our approach allows for the determination of the simplex packing widths in dimension six for all $k\le8$.

In addition, our approach determines all $k$-cardinality simplex-multisets that allow for optimal packings. For each such multiset, an explicit optimal packing is generated.
We present and discuss the number of optimal multisets as a function of~$k$.

In Section~\ref{sec: MILP1} we describe how we model the determination of the simplex packing width $s_k \of{P^4(1), \omega_0}$ as a mixed integer linear program (MILP), along with the presentation and discussion of the computer implementation and the computational results.
In Section~\ref{sec: SDP} we present an alternative modelling as a quadratically constrained quadratic program (QCQP) which is then relaxed to a semidefinite program (SDP). In Section~\ref{sec: MILP2} we extend the MILP approach to the next higher dimension and compute the $k$-simplex packing width of the six-dimensional open prism.

\section{MILP for $k$-Simplex Packing Width of $P^4(1)$}\label{sec: MILP1}

\subsection*{Modelling}

Whereas algebraic geometry is a crucial tool to calculate the ball packing widths, the main tool to calculate the simplex packing widths is the following theorem. 
\begin{Thm}[Packing Theorem~\cite{Maley2000}]\label{thm: packing theorem}
    Let $V$ be an open, connected subset of $\R^n$ with $H_1(V,\Z) = 0$. Then 
    \begin{equation*}
        s_k \of{\T^n \times V, \omega_0} = \operatorname{sup} \offf{r \st{ 
        \begin{array}{l}
            \exists A_1, \ldots, A_k \in \GL{n} \; \exists t_1, \ldots, t_k \in \R^n : \\
            A_i \of{\standardsimplex^n(r)} + t_i \subseteq V \quad \forall i \in \offf{1,\ldots,k}\\
            \of{A_i \of{\standardsimplex^n(r)} + t_i} \cap \of{A_j \of{\standardsimplex^n(r)} + t_j} = \emptyset \\
            \forall 1 \leq i < j \leq k
        \end{array} 
        }}.
    \end{equation*}
\end{Thm}
Here, $\GL{n}$ denotes the set of all matrices in $\Z^{n\times n}$ that are invertible over~$\Z$, together with matrix multiplication as the group operation. The meaning of Theorem~\ref{thm: packing theorem} is that for symplectic manifolds of the form $\T^n \times V$, where $V$ is an open connected subset of $\R^n$ with first homology equal to zero, we can compute the $k$-simplex packing width $s_k \of{\T^n \times V, \omega_0}$ by computing an optimal packing of $V$ by copies of $\standardsimplex^n(r)$ under integral affine transformations while maximizing the sidelength $r$. This not only reduces the dimension of the problem space from $2n$ to $n$ but also converts the calculation of the $k$-simplex packing width into a classical combinatoral packing problem. By applying Theorem~\ref{thm: packing theorem} to $P^4(1)=\T^2\times\standardsimplex^2(1)$, we get the following  formulation.
\begin{Prb}[Outer Optimization Problem - Combinatorial Formulation]\label{prb: outer2}
    Given $k \in \N$, determine the minimum side length $s$ such that there exist matrices $A_1, \ldots, A_k \in \GL{2}$ and vectors $t_1, \ldots, t_k \in \R^2$ satisfying the \emph{containment condition}
\[A_i \of{\standardsimplex^2(1)} + t_i \subseteq \standardsimplex^2(s)\quad\mbox{ for all }1\le i\le k\]
and the \emph{disjointness condition}
\[\of{A_i \of{\standardsimplex^2(1)} + t_i} \cap \of{A_j \of{\standardsimplex^2(1)} + t_j} = \emptyset\quad\mbox{ for all }1 \leq i < j \leq k.\]
\end{Prb}

Notice that previously we were maximizing the size of the packing objects but now we are minimizing the size of the packing container. The reason for this is that the latter is easier to model.
Let $s_k^{\standardsimplex}$ denote the minimum side length~$s$ from Problem~\ref{prb: outer2}. Then we have the relation 
\begin{equation*}
    s_k \of{P^4(1), \omega_0} = \frac{1}{s_k^{\standardsimplex}}.
\end{equation*} 
By taking the reciprocal value of $s_k \of{P^4(1), \omega_0}$ in Table~\ref{tab: packing-widths}, we obtain an upper bound $\overline{s_k^{\standardsimplex}}$ on $s_k^{\standardsimplex}$ for $k=1,\ldots,20$. For $k \leq 12$ and $k = 16$ we even have equality $\overline{s_k^{\standardsimplex}} = s_k^{\standardsimplex}$. We only need to consider triangles whose shapes are admissible for that upper bound. For a given $k$, we call the list of admissible triangles the shapelist $\mathcal{S}^{\standardsimplex}_k$. Now we can determine $s_k^{\standardsimplex}$ by computing an optimal packing for every $k$-cardinality multisubset of the shapelist. The number of $k$-cardinality multisubsets from the shapelist $\of{\vector{\abs{\mathcal{S}_k^{\standardsimplex}}}{k}}$ according to $k$ is shown in Table~\ref{tab: multisubsets}. 

\begin{table}[ht]
    \centering
    \begin{small}
    \begin{tabular}{|R|R|R|}                                 
        \hline
        \rule{0pt}{15pt} k & \multicolumn{1}{C}{\abs{\mathcal{S}_k^{\standardsimplex}}} & \multicolumn{1}{|C|}{\of{\vector{\abs{\mathcal{S}_k^{\standardsimplex}}}{k}}} \\[7pt]
        \hline
        1 & 1 & 1\\
        2 & 8 & 36\\
        3 & 8 & 120\\
        4 & 8 & 330\\
        5 & 8 & 792\\
        6 & 8 & 1\,716\\
        7 & 20 & 657\,800\\
        8 & 20 & 2\,220\,075\\
        9 & 20 & 6\,906\,900\\
        10 & 20 & 20\,030\,010\\
        11 & 20 & 54\,627\,300\\
        12 & 20 & 141\,120\,525\\
        13 & 20 & 347\,373\,600\\
        14 & 20 & 818\,809\,200\\
        15 & 32 & 511\,738\,760\,544\\
        16 & 32 & 1\,503\,232\,609\,098\\
        17 & 32 & 4\,244\,421\,484\,512\\
        18 & 32 & 11\,554\,258\,485\,616\\
        19 & 32 & 30\,405\,943\,383\,200\\
        20 & 32 & 77\,535\,155\,627\,160\\
        \hline
    \end{tabular}
    \end{small}
    \caption{Number of $k$-cardinality multisubsets of the shapelists}
    \label{tab: multisubsets}
\end{table}

For $k = 10$ one already needs to solve more than twenty million subproblems. This clearly indicates that complete enumeration is out of the question, especially in view of the fact that only exponential time algorithms for the computation of the optimal packings are known. Rather, like~\cite{Maley2000} we implement a branch-and-bound approach. Each node of the branch-and-bound search tree represents a multisubset. The level of the node is equal to the cardinality of the multisubset. For example, the root on level zero corresponds to the empty multisubset and the children on level one correspond to all subsets of cardinality one.

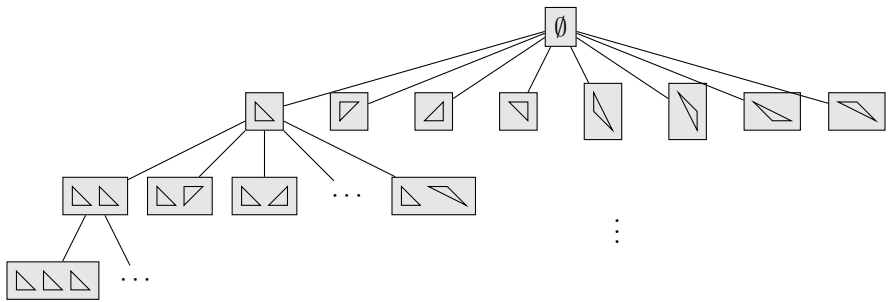
\begin{figure}[ht]
    \centering
    \begin{tikzpicture}[scale = 0.75, standard/.style = {draw, fill=gray!20}, dots/.style = {}]
        \node[standard] {$\emptyset$}
            child {node[standard] {$\triang{0}{1}{1}{0}$}
                child {node[standard] {$\triang{0}{1}{1}{0} \; \triang{0}{1}{1}{0}$}
                    child {node[standard] {$\triang{0}{1}{1}{0} \; \triang{0}{1}{1}{0} \; \triang{0}{1}{1}{0}$}}
                    child {node[dots] {$\cdots$}}}
                child {node[standard] {$\triang{0}{1}{1}{0} \; \triang{0}{1}{1}{1}$}}
                child {node[standard] {$\triang{0}{1}{1}{0} \; \triang{1}{0}{1}{1}$}}
                child {node[dots] {$\cdots$}}
                child {node[standard] {$\triang{0}{1}{1}{0} \; \triang{1}{0}{2}{-1}$}}}
            child {node[standard] {$\triang{0}{1}{1}{1}$}}
            child {node[standard] {$\triang{1}{0}{1}{1}$}} 
            child {node[standard] {$\triang{1}{-1}{1}{0}$}} 
            child {node[standard] {$\triang{0}{1}{1}{-1}$}} 
            child {node[standard] {$\triang{1}{-2}{1}{-1}$}} 
            child {node[standard] {$\triang{1}{-1}{2}{-1}$}} 
            child {node[standard] {$\triang{1}{0}{2}{-1}$}};
        \draw node[above] at (1,-4) {$\vdots$};
    \end{tikzpicture}
    \caption{Branch-and-bound search tree for the $3$-triangle packing}
    \label{fig: search tree}
\end{figure}

The (incomplete) search tree for $k = 3$ is shown in Figure~\ref{fig: search tree}.

We start the search with the global upper bound obtained from Table~\ref{tab: packing-widths}. Whenever the computation for a node produces an optimum packing that exceeds the global upper bound, we can discard 
the subtree rooted at this node
from further investigation: the node (together with its subtree) is \emph{fathomed}. Whenever the computation for a node at level $k$ produces an optimum packing that improves the global upper bound, we update the global upper bound and memorize the packing as the incumbent solution. At termination the incumbent solution is an optimum packing with value of the upper bound.

Now we address how to 
compute an optimal packing at each node. We call this the inner optimization problem.
\begin{Prb}[Inner Optimization Problem]\label{prb: inner1}
    Given $T_1, \ldots, T_m \in \mathcal{S}_k^{\standardsimplex}$, determine the minimum side length $s$ such that there exist vectors $t_1, \ldots, t_m \in \R^2$ satisfying 
   \begin{align*}
        &T_i + t_i \subseteq \standardsimplex^2(s) &\mbox{for all } &i \in \offf{1,\ldots,m} &\text{(containment condition)},\\
        &\of{T_i + t_i} \cap \of{T_j + t_j} = \emptyset &\mbox{for all } &1 \leq i < j \leq m &\text{(disjointness condition)}.
    \end{align*}
\end{Prb}
We 
formulate Problem~\ref{prb: inner1} as a MILP. Every triangle from the shapelist is given in the form 
\begin{equation*}
    T_i = \conv{0}{0}{a_i}{b_i}{c_i}{d_i},
\end{equation*}
where $a_i,b_i,c_i,d_i$ are integer constants.
Instead of working with the open sets $T_i$ and $\standardsimplex^2(s)$, we will consider their closures. This does not make any difference for the containment condition but is easier to model. Let $t_i = \vector{x_i}{y_i}$ denote the translation vector, then the closure of the translated triangle $T_i + t_i$ is given by
\begin{equation*}
    \overline{T_i + t_i} = \conv{x_i}{y_i}{x_i + a_i}{y_i + b_i}{x_i + c_i}{y_i + d_i}
\end{equation*}
and the closed standard simplex of side length $s$ is given by
\begin{equation*}
    \overline{\standardsimplex^2(s)} = \offf{\left. \vector{x}{y} \in \R^2 \;\right|\; x \geq 0,\; y \geq 0,\; x + y \leq s}.
\end{equation*}
Because of convexity, it suffices to check the containment condition for the three vertices only. Thus, we obtain a total of nine inequalities for every $i \in \offf{1,\ldots,m}$:
\begin{align*}
    x_i       &\geq 0,  & y_i       &\geq 0,    & x_i + y_i                 &\leq s, \\ 
    x_i + a_i &\geq 0,  & y_i + b_i &\geq 0,    & (x_i + a_i) + (y_i + b_i) &\leq s, \\ 
    x_i + c_i &\geq 0,  & y_i + d_i &\geq 0,    & (x_i + c_i) + (y_i + d_i) &\leq s.    
\end{align*}
By putting the constants to the right hand sides and taking extrema in every column, we can reduce these nine inequalities 
%to three inequalities 
to only three
for every $i \in \of{1,\ldots,m}$:

\begin{align*}
    x_i             &\geq \mmax{0, -a_i, -c_i}              &=:\;& K_i, \\ 
    y_i             &\geq \mmax{0, -b_i, -d_i}              &=:\;&K'_i,\\ 
    x_i + y_i -s    &\leq \mmin{0, -a_i - b_i, -c_i - d_i}  &=:\;&K''_i.
\end{align*}%
Thus, we get a total of $3m$ inequalities that we call the containment constraints.

Next, we model the disjointness condition. We observe that the equation $\of{T_i + t_i} \cap \of{T_j + t_j} = \emptyset$
is equivalent to $t_j - t_i \notin \of{T_i \ominus T_j}$. Here, $T_i \ominus T_j$ denotes the Minkowski difference of $T_i$ and $T_j$ that is defined as 
\begin{equation*}
   T_i \ominus T_j = \offf{v_i - v_j \st v_i \in T_i, v_j \in T_j}. 
\end{equation*}
In our setting, the Minkowski difference of two triangles $T_i$ and $T_j$ is of the form 
\begin{equation*}
    T_i \ominus T_j = \offf{ \left. \vector{x}{y} \in \R^2 \;\right|\; \alpha^{ij}_f x + \beta^{ij}_f y < \gamma^{ij}_f \quad \forall f\in\offf{1,\ldots,6}}, 
\end{equation*} 
where $\alpha^{ij}_f, \beta^{ij}_f\gamma^{ij}_f \in \Z$ are integer constants. Working with this representation, the difference vector $t_j - t_i$ is not contained in the Minkowski difference $T_i \ominus T_j$ if and only if at least one of the six inequalities $\alpha^{ij}_f x + \beta^{ij}_f y < \gamma^{ij}_f$ is violated. To model this condition, we introduce a binary variable $z^{ij}_f \in \offf{0,1}$ for every $f \in \offf{1,\ldots,6}$ with the following meaning:
\begin{equation*}
    z^{ij}_f = 1 \Rightarrow \alpha^{ij}_f (x_j - x_i) + \beta^{ij}_f (y_j - y_i) \geq \gamma^{ij}_f.
\end{equation*}
Now one could by brute force enumerate all $\mathcal{O}\of{6^{\binom{m}{2}}}$ possible $0/1$-assignments of the binary variables $z^{ij}_f$ and solve a linear program for each. This strategy was pursued by Maley et~al.~\cite{Maley2000}. Instead, we model the implication using a Big-$M$-formulation, where the parameter $M$ has to be chosen sufficiently large:
\begin{equation*}
    \alpha^{ij}_f (x_j - x_i) + \beta^{ij}_f (y_j - y_i) \geq  \gamma^{ij}_f - M \of{ 1 - z^{ij}_f }.
\end{equation*}
In our setting we can choose $M$ as 
\begin{equation*}
    M = \of{\abs{\alpha^{ij}_f} + \abs{\beta^{ij}_f}} \hat{s} + \gamma^{ij}_f,
\end{equation*}
where $\hat{s}$ is the current global upper bound on $s$ in the branch-and-bound algorithm when encountering the inner optimization problem.

We want at least one of the six binary variables $z^{ij}_f$ to take the value one. Equivalently we can say that their sum should be greater or equal to one. Altogether, the disjointness condition is equivalent to
\begin{align*}
    &z^{ij}_f \in \offf{0,1},\\ 
    &z^{ij}_1 + \cdots + z^{ij}_6 \geq 1,\\ 
    &\alpha^{ij}_f (x_j - x_i) + \beta^{ij}_f (y_j - y_i) \geq  \gamma^{ij}_f - M \of{ 1 - z^{ij}_f }
\end{align*}
for all $f \in \offf{1,\ldots,6}$. For every pair of triangles $T_i,T_j$ with $1 \leq i < j \leq m$, we have these seven inequalities, so in total there are $7\binom{m}{2}$ inequalities that we call the disjointness constraints. 
By merging the containment constraints and the disjointness constraints, we obtain an MILP with $1 + 2m$ continuous variables, $6\binom{m}{2}$ (binary) integer variables and $3m + 7\binom{m}{2}$ constraints. 
\begin{Prb}[MILP Formulation of the Inner Optimization Problem~\ref{prb: inner1}]\label{prb: inner2}
    \begin{align*} 
        \min \quad & s & & \\
        & x_i \geq K_i & \forall & i \in \offf{1,\ldots,m} \\ 
        & y_i \geq K'_i & \forall & i \in \offf{1,\ldots,m} \\ 
        & x_i + y_i - s \leq K''_i & \forall & i \in \offf{1,\ldots,m} \\
        & z^{ij}_1 + \cdots + z^{ij}_6 \geq 1 & \forall & 1 \leq i < j \leq m \\
        & \alpha^{ij}_f (x_j - x_i) + \beta^{ij}_f (y_j - y_i) \geq  \gamma^{ij}_f - M \of{ 1 - z^{ij}_f } & \forall & 1 \leq i < j \leq m \\
        & &\forall & f \in \offf{1,\ldots,6} \\
        & s \in \R, \quad x,y \in \R^m, \quad z \in \offf{0, 1}^{6\binom{m}{2}} 
    \end{align*} 
\end{Prb}

\subsection*{Implementation}

We keep a data base of multisets with known lower bounds or optimal solutions, initialized with the empty set. Then we run our program with input $k$ for $k=1,2,\ldots$. Via the data base the run for $k$ can profit from data base queries for given multisets of cardinalities up to $k-1$ rather than solving the inner optimization problem.

The program maintains a list of \emph{blocking multisets}, i.e., multisets for which a lower bound exceeding the global upper bound is known. 

The backbone of the program is a queue of branch-and-bound nodes. Each of these nodes represents an inner optimization instance. As long as the queue is not empty, we extract its front element. If the lower bound of the inner optimization instance exceeds the global upper bound, the node is fathomed. Otherwise we check the inner optimization instance for blocking multisubsets. If we detect a blocking subset, the node is fathomed. Otherwise, we check for previously known bounds. If a lower bound is found that exceeds the global upper bound, the inner optimization instance is recorded as a new blocking structure and the node is fathomed. In any other case, we run the inner optimization solver
and update the database of known solutions and lower bounds. 

If the solution computed by the inner optimization solver exceeds the global upper bound, the inner optimization instance is recorded as a new blocking structure and the node is fathomed. If the optimal value of the solution is smaller or equal to the global upper bound, we distinguish between two cases. 

In the first case, the number of simplices in the inner optimization instance is equal to the number $k$ of simplices in the outer optimization problem, which means the node is a leaf of the branch-and-bound tree. If the optimal value of the solution is equal to the current global upper bound, we draw a picture of the solution and append it to a \LaTeX-file. This \LaTeX-file already contains pictures of all previous solutions with the same optimal value. If the optimal value of the solution is strictly smaller than the global upper bound, we update the latter and replace the old \LaTeX-file by a new \LaTeX-file containing a picture of the solution.

In the second case, the number of simplices in the inner optimization instance is smaller than the number $k$ of simplices in the outer optimization problem, which means the node is not a leaf. If so, we add all feasible extensions of the inner optimization instance to the queue. 

\begin{figure}[!ht]
    \centering
    \begin{tikzpicture}[S/.style = draw, rounded corners, align=center, scale = 0.67]
    \linespread{1}
    \def\l{0}
    \def\m{6}
    \def\r{12}
    \def\h{-3}
    
        \node at (\l,0) [S] (l0) {Extract the front\\ element of the queue.};
        \node at (\r,2*\h) [S] (r) {Record new\\ blocking\\ configuration.};
        \node at (\r,\h) [S] (f) {Fathom.};
        \draw[->,thick] (r) -- (f);
        
        \node at (\l,\h) [S] (l1) {Check for\\ blocking\\ configurations.};
        \draw[->,thick] (l0) -- (l1);
        
        \node at (\l,2*\h) [S] (l2) {Check for\\ known\\ lower bounds.};
        \draw[->,thick] (l1) -- (l2) node[midway, left] {\scriptsize{negative}};
        \draw[->,thick] (l1) -- (f) node[midway, below] {\scriptsize{positive}};
        
        \node at (\l,3*\h) [S] (l3) {Run the inner\\ optimization solver\\ and append the result\\ to the bounds file.};
        \draw[->,thick] (l2) -- (l3) node[midway, left] {\scriptsize{negative}};
        \node at (\m,2*\h) [S] (m3) {Does the lower\\ bound exceed the\\ global upper bound ?};
        \draw[->,thick] (l2) -- (m3) node[midway, below] {\scriptsize{positive}};
        \draw[->,thick] (m3) -- (r) node[midway, below] {\scriptsize{yes}};
        \draw[->,thick] (m3) -- ++(0,\h) node[midway, left] {\scriptsize{no}} -- (l3);
        
        \node at (\l,4*\h) [S] (l4) {Is the problem\\ feasible?};
        \draw[->,thick] (l3) -- (l4);
        
        \node at (\l,5*\h) [S] (l5) {Does the solution\\ exceed the global\\ upper bound ?};
        \draw[->,thick] (l4) -- (l5) node[midway, left] {\scriptsize{yes}};
        \draw[->,thick] (l4) -- ++ (\r,0) node[midway, below] {\scriptsize{no}} -- (r);

        \node at (\l,6*\h) [S] (l6) {Is the current\\ number of simplices\\ equal to $k$?};
        \draw[->,thick] (l5) -- (l6) node[midway, left] {\scriptsize{no}};        
        \draw[->,thick] (l5) -- ++(\r,0) node[midway, below] {\scriptsize{yes}} -- (r);
        
        \node at (\l,7*\h) [S] (l7) {Add all feasible\\ extensions \\to the queue.};
        \draw[->,thick] (l6) -- (l7) node[midway, left] {\scriptsize{no}};        
        \node at (\m,6*\h) [S] (m7) {Update the global\\ upper bound and \\ write a picture of the \\ solution to a \LaTeX-file.};
        \draw[->,thick] (l6) -- (m7) node[midway, below] {\scriptsize{yes}};
        
    \end{tikzpicture}
    \caption{Flowchart of our program}
    \label{fig: workflow}
\end{figure}
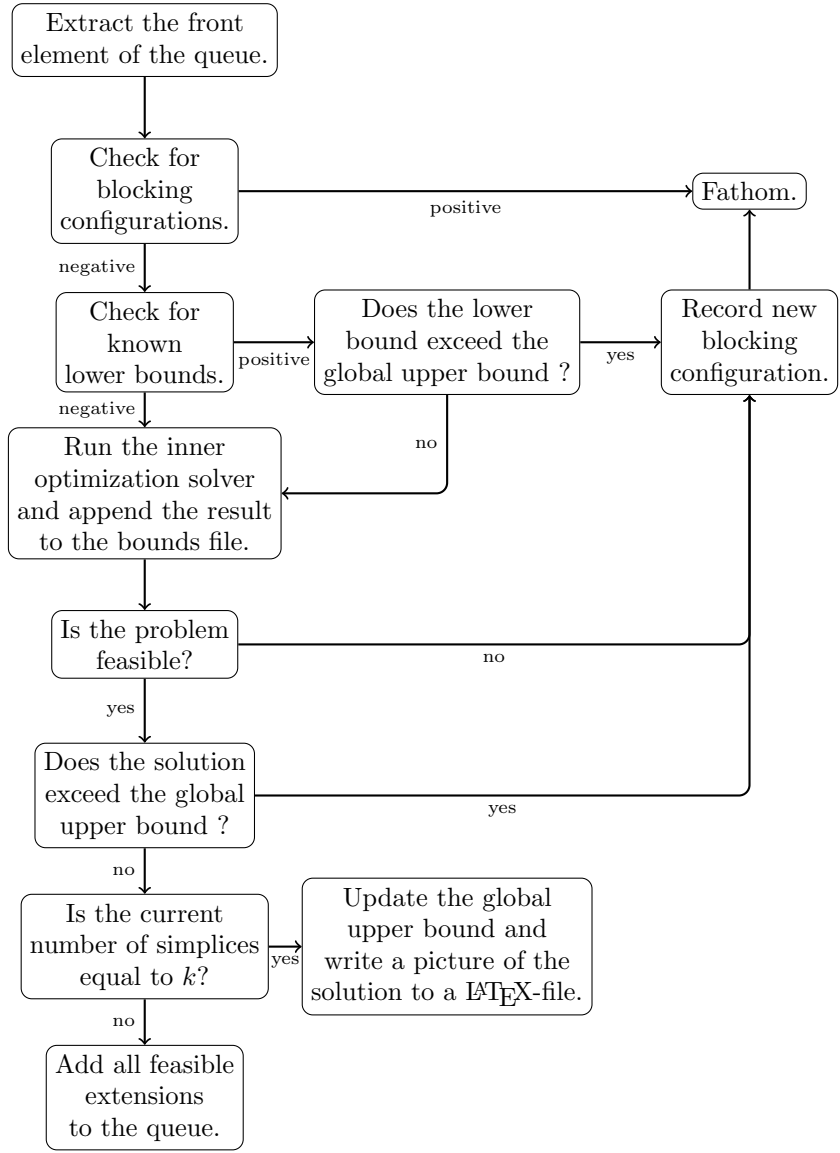

All described steps are summarized in Figure~\ref{fig: workflow}. 
The computer code consists of
1841 lines
written in the programming language C,
406 
of which
($\approx 22\%$)
concern the branch-and-bound framework.
The code for the inner optimization procedure 
has 471 lines ($\approx26\%$). 
Reading and rewriting the bounds file and the \LaTeX-file %it 
takes 204 lines ($\approx11\%$) and 162 lines ($\approx9\%$), respectively. 

\subsection*{Results}

The computational experiments were carried out under the Debian 10 operating system on two Intel Xeon E5-2690v2 CPUs with 3.00GHz and 10 cores each. For solving the inner optimization problem, we call the GUROBI Optimization Software Version 9.0.3. All computation times are wall clock times.

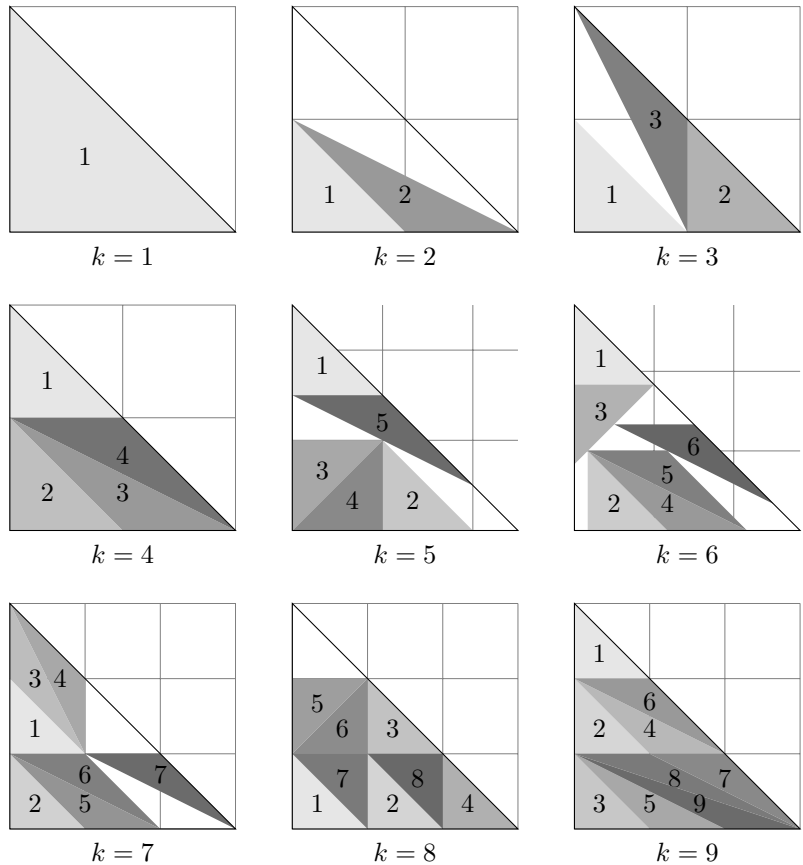
\begin{figure}[ht]
    \centering
    \begin{minipage}{0.3\textwidth}
        \centering
        \begin{tikzpicture}[scale=3] 
            \draw [help lines] (0.0,0.0) grid (1.000000,1.000000); 
            \fill[black!10] (0.000000,0.000000)--(0.000000,1.000000)--(1.000000,0.000000)--cycle;
            \node at (0.333333,0.333333) {1}; 
            \draw (0.0,0.0)--(1.000000,0.0)--(0.0,1.000000)--cycle; 
        \end{tikzpicture}
        \\$k=1$
    \end{minipage}
    \begin{minipage}{0.3\textwidth}
        \centering
        \begin{tikzpicture}[scale=1.5] 
            \draw [help lines] (0.0,0.0) grid (2.000000,2.000000); 
            \fill[black!10] (0.000000,0.000000)--(0.000000,1.000000)--(1.000000,0.000000)--cycle;
            \node at (0.333333,0.333333) {1}; 
            \fill[black!40] (0.000000,1.000000)--(1.000000,0.000000)--(2.000000,0.000000)--cycle;
            \node at (1.000000,0.333333) {2}; 
            \draw (0.0,0.0)--(2.000000,0.0)--(0.0,2.000000)--cycle; 
        \end{tikzpicture}
        \\$k=2$
    \end{minipage}
    \begin{minipage}{0.3\textwidth}
        \centering
        \begin{tikzpicture}[scale=1.5] 
            \draw [help lines] (0.0,0.0) grid (2.000000,2.000000); 
            \fill[black!10] (0.000000,0.000000)--(0.000000,1.000000)--(1.000000,0.000000)--cycle;
            \node at (0.333333,0.333333) {1}; 
            \fill[black!30] (1.000000,0.000000)--(1.000000,1.000000)--(2.000000,0.000000)--cycle;
            \node at (1.333333,0.333333) {2}; 
            \fill[black!50] (0.000000,2.000000)--(1.000000,0.000000)--(1.000000,1.000000)--cycle;
            \node at (0.72,1.000000) {3}; 
            \draw (0.0,0.0)--(2.000000,0.0)--(0.0,2.000000)--cycle; 
        \end{tikzpicture}
        \\$k=3$
    \end{minipage}
    \\[0.5cm]
    \begin{minipage}{0.3\textwidth}
        \centering
        \begin{tikzpicture}[scale=1.5] 
            \draw [help lines] (0.0,0.0) grid (2.000000,2.000000); 
            \fill[black!10] (0.000000,1.000000)--(0.000000,2.000000)--(1.000000,1.000000)--cycle;
            \node at (0.333333,1.333333) {1}; 
            \fill[black!25] (0.000000,0.000000)--(0.000000,1.000000)--(1.000000,0.000000)--cycle;
            \node at (0.333333,0.333333) {2}; 
            \fill[black!40] (0.000000,1.000000)--(1.000000,0.000000)--(2.000000,0.000000)--cycle;
            \node at (1.000000,0.333333) {3}; 
            \fill[black!55] (0.000000,1.000000)--(1.000000,1.000000)--(2.000000,0.000000)--cycle;
            \node at (1.000000,0.666667) {4}; 
            \draw (0.0,0.0)--(2.000000,0.0)--(0.0,2.000000)--cycle; 
        \end{tikzpicture}
        \\$k=4$
    \end{minipage}
    \begin{minipage}{0.3\textwidth}
        \centering
        \begin{tikzpicture}[scale=1.2] 
            \draw [help lines] (0.0,0.0) grid (2.500000,2.500000); 
            \fill[black!10] (0.000000,1.500000)--(0.000000,2.500000)--(1.000000,1.500000)--cycle;
            \node at (0.333333,1.833333) {1}; 
            \fill[black!22] (1.000000,0.000000)--(1.000000,1.000000)--(2.000000,0.000000)--cycle;
            \node at (1.333333,0.333333) {2}; 
            \fill[black!34] (0.000000,0.000000)--(0.000000,1.000000)--(1.000000,1.000000)--cycle;
            \node at (0.333333,0.666667) {3}; 
            \fill[black!46] (0.000000,0.000000)--(1.000000,0.000000)--(1.000000,1.000000)--cycle;
            \node at (0.666667,0.333333) {4}; 
            \fill[black!58] (0.000000,1.500000)--(1.000000,1.500000)--(2.000000,0.500000)--cycle;
            \node at (1.000000,1.2) {5}; 
            \draw (0.0,0.0)--(2.500000,0.0)--(0.0,2.500000)--cycle; 
        \end{tikzpicture}
        \\$k=5$
    \end{minipage}
    \begin{minipage}{0.3\textwidth}
        \centering
        \begin{tikzpicture}[scale=1.05882353] 
            \draw [help lines] (0.0,0.0) grid (2.833333,2.833333); 
            \fill[black!10] (0.000000,1.833333)--(0.000000,2.833333)--(1.000000,1.833333)--cycle;
            \node at (0.333333,2.166667) {1}; 
            \fill[black!20] (0.166667,0.000000)--(0.166667,1.000000)--(1.166667,0.000000)--cycle;
            \node at (0.500000,0.333333) {2}; 
            \fill[black!30] (0.000000,0.833333)--(0.000000,1.833333)--(1.000000,1.833333)--cycle;
            \node at (0.333333,1.500000) {3}; 
            \fill[black!40] (0.166667,1.000000)--(1.166667,0.000000)--(2.166667,0.000000)--cycle;
            \node at (1.166667,0.333333) {4}; 
            \fill[black!50] (0.166667,1.000000)--(1.166667,1.000000)--(2.166667,0.000000)--cycle;
            \node at (1.166667,0.7) {5}; 
            \fill[black!60] (0.500000,1.333333)--(1.500000,1.333333)--(2.500000,0.333333)--cycle;
            \node at (1.500000,1.05) {6}; 
            \draw (0.0,0.0)--(2.833333,0.0)--(0.0,2.833333)--cycle; 
        \end{tikzpicture}
        \\$k=6$
    \end{minipage}
    \\[0.5cm]
    \begin{minipage}{0.3\textwidth}
    \centering
    \begin{tikzpicture}[scale=1.0] 
        \draw [help lines] (0.0,0.0) grid (3.000000,3.000000); 
        \fill[black!10] (0.000000,1.000000)--(0.000000,2.000000)--(1.000000,1.000000)--cycle;
        \node at (0.333333,1.333333) {1}; 
        \fill[black!18] (0.000000,0.000000)--(0.000000,1.000000)--(1.000000,0.000000)--cycle;
        \node at (0.333333,0.333333) {2}; 
        \fill[black!26] (0.000000,2.000000)--(0.000000,3.000000)--(1.000000,1.000000)--cycle;
        \node at (0.333333,2.000000) {3}; 
        \fill[black!34] (0.000000,3.000000)--(1.000000,1.000000)--(1.000000,2.000000)--cycle;
        \node at (0.666667,2.000000) {4}; 
        \fill[black!42] (0.000000,1.000000)--(1.000000,0.000000)--(2.000000,0.000000)--cycle;
        \node at (1.000000,0.333333) {5}; 
        \fill[black!50] (0.000000,1.000000)--(1.000000,1.000000)--(2.000000,0.000000)--cycle;
        \node at (1.000000,0.72) {6}; 
        \fill[black!58] (1.000000,1.000000)--(2.000000,1.000000)--(3.000000,0.000000)--cycle;
        \node at (2.000000,0.72) {7}; 
        \draw (0.0,0.0)--(3.000000,0.0)--(0.0,3.000000)--cycle; 
    \end{tikzpicture}
    \\$k=7$
    \end{minipage}
        \begin{minipage}{0.3\textwidth}
        \centering
        \begin{tikzpicture}[scale=1] 
            \draw [help lines] (0.0,0.0) grid (3.000000,3.000000); 
                \fill[black!10] (0.000000,0.000000)--(0.000000,1.000000)--(1.000000,0.000000)--cycle;
                \node at (0.333333,0.333333) {1}; 
                \fill[black!17] (1.000000,0.000000)--(1.000000,1.000000)--(2.000000,0.000000)--cycle;
                \node at (1.333333,0.333333) {2}; 
                \fill[black!24] (1.000000,1.000000)--(1.000000,2.000000)--(2.000000,1.000000)--cycle;
                \node at (1.333333,1.333333) {3}; 
                \fill[black!31] (2.000000,0.000000)--(2.000000,1.000000)--(3.000000,0.000000)--cycle;
                \node at (2.333333,0.333333) {4}; 
                \fill[black!38] (0.000000,1.000000)--(0.000000,2.000000)--(1.000000,2.000000)--cycle;
                \node at (0.333333,1.666667) {5}; 
                \fill[black!45] (0.000000,1.000000)--(1.000000,1.000000)--(1.000000,2.000000)--cycle;
                \node at (0.666667,1.333333) {6}; 
                \fill[black!52] (0.000000,1.000000)--(1.000000,0.000000)--(1.000000,1.000000)--cycle;
                \node at (0.666667,0.666667) {7}; 
                \fill[black!59] (1.000000,1.000000)--(2.000000,0.000000)--(2.000000,1.000000)--cycle;
                \node at (1.666667,0.666667) {8}; 
                \draw (0.0,0.0)--(3.000000,0.0)--(0.0,3.000000)--cycle; 
        \end{tikzpicture}
        \\$k=8$
    \end{minipage}
    \begin{minipage}{0.3\textwidth}
        \centering
        \begin{tikzpicture}[scale=1] 
            \draw [help lines] (0.0,0.0) grid (3.000000,3.000000); 
            \fill[black!10] (0.000000,2.000000)--(0.000000,3.000000)--(1.000000,2.000000)--cycle;
            \node at (0.333333,2.333333) {1}; 
            \fill[black!16] (0.000000,1.000000)--(0.000000,2.000000)--(1.000000,1.000000)--cycle;
            \node at (0.333333,1.333333) {2}; 
            \fill[black!22] (0.000000,0.000000)--(0.000000,1.000000)--(1.000000,0.000000)--cycle;
            \node at (0.333333,0.333333) {3}; 
            \fill[black!28] (0.000000,2.000000)--(1.000000,1.000000)--(2.000000,1.000000)--cycle;
            \node at (1.000000,1.333333) {4}; 
            \fill[black!34] (0.000000,1.000000)--(1.000000,0.000000)--(2.000000,0.000000)--cycle;
            \node at (1.000000,0.333333) {5}; 
            \fill[black!40] (0.000000,2.000000)--(1.000000,2.000000)--(2.000000,1.000000)--cycle;
            \node at (1.000000,1.7) {6}; 
            \fill[black!46] (1.000000,1.000000)--(2.000000,1.000000)--(3.000000,0.000000)--cycle;
            \node at (2.000000,0.666667) {7}; 
            \fill[black!52] (0.000000,1.000000)--(1.000000,1.000000)--(3.000000,0.000000)--cycle;
            \node at (1.333333,0.666667) {8}; 
            \fill[black!58] (0.000000,1.000000)--(2.000000,0.000000)--(3.000000,0.000000)--cycle;
            \node at (1.666667,0.333333) {9}; 
            \draw (0.0,0.0)--(3.000000,0.0)--(0.0,3.000000)--cycle; 
        \end{tikzpicture}
        \\$k=9$
    \end{minipage}
    \caption{Optimal $k$-triangle packings for $k = 1,\ldots,9$}
    \label{fig: packings1}
\end{figure}

\begin{figure}[ht]
    \centering   
    \begin{minipage}{0.45\textwidth}
        \centering
        %\begin{tikzpicture}[scale=1.2] 
        %\begin{tikzpicture}[scale=1.44] 
        %\begin{tikzpicture}[scale=1.56] 
        \begin{tikzpicture}[scale=1.5] 
            \draw [help lines] (0.0,0.0) grid (3.333332,3.333332); 
            \fill[black!10] (0.000000,0.000000)--(0.000000,1.000000)--(1.000000,0.000000)--cycle;
            \node at (0.333333,0.333333) {1}; 
            \fill[black!16] (1.333333,0.000000)--(1.333333,1.000000)--(2.333333,0.000000)--cycle;
            \node at (1.666667,0.333333) {2}; 
            \fill[black!22] (1.333333,1.000000)--(1.333333,2.000000)--(2.333333,1.000000)--cycle;
            \node at (1.666667,1.333333) {3}; 
            \fill[black!28] (0.000000,1.000000)--(0.000000,2.000000)--(1.000000,0.000000)--cycle;
            \node at (0.333333,1.000000) {4}; 
            \fill[black!34] (0.333333,3.000000)--(1.333333,1.000000)--(1.333333,2.000000)--cycle;
            \node at (1.000000,2.000000) {5}; 
            \fill[black!40] (1.333333,1.000000)--(2.333333,0.000000)--(3.333333,0.000000)--cycle;
            \node at (2.333333,0.333333) {6}; 
            \fill[black!46] (1.333333,1.000000)--(2.333333,1.000000)--(3.333333,0.000000)--cycle;
            \node at (2.333333,0.666667) {7}; 
            \fill[black!52] (0.000000,2.000000)--(0.000000,3.000000)--(1.000000,0.000000)--cycle;
            \node at (0.333333,1.666667) {8}; 
            \fill[black!58] (0.000000,3.000000)--(1.000000,0.000000)--(1.000000,1.000000)--cycle;
            \node at (0.666667,1.333333) {9}; 
            \fill[black!64] (0.333333,3.000000)--(1.333333,0.000000)--(1.333333,1.000000)--cycle;
            \node at (1.000000,1.333333) {10}; 
            \draw (0.0,0.0)--(3.333332,0.0)--(0.0,3.333332)--cycle; 
        \end{tikzpicture}
        \\$k=10$
    \end{minipage}
    \begin{minipage}{0.45\textwidth}
        \centering
        %\begin{tikzpicture}[scale=1.14285714] 
        %\begin{tikzpicture}[scale=1.37145017] 
        %\begin{tikzpicture}[scale=1.48573768] 
        \begin{tikzpicture}[scale=1.42859392] 
            \draw [help lines] (0.0,0.0) grid (3.500000,3.500000); 
            \fill[black!10] (0.000000,0.000000)--(0.000000,1.000000)--(1.000000,0.000000)--cycle;
            \node at (0.333333,0.333333) {1}; 
            \fill[black!15] (1.000000,1.500000)--(1.000000,2.500000)--(2.000000,1.500000)--cycle;
            \node at (1.333333,1.833333) {2}; 
            \fill[black!20] (0.000000,2.500000)--(0.000000,3.500000)--(1.000000,2.500000)--cycle;
            \node at (0.333333,2.833333) {3}; 
            \fill[black!25] (0.000000,1.500000)--(0.000000,2.500000)--(1.000000,2.500000)--cycle;
            \node at (0.333333,2.166667) {4}; 
            \fill[black!30] (0.000000,1.500000)--(1.000000,1.500000)--(1.000000,2.500000)--cycle;
            \node at (0.666667,1.833333) {5}; 
            \fill[black!35] (0.000000,1.000000)--(1.000000,0.000000)--(2.000000,0.000000)--cycle;
            \node at (1.000000,0.333333) {6}; 
            \fill[black!40] (1.000000,1.500000)--(2.000000,1.500000)--(3.000000,0.500000)--cycle;
            \node at (2.000000,1.166667) {7}; 
            \fill[black!45] (0.000000,1.000000)--(1.000000,1.000000)--(3.000000,0.000000)--cycle;
            \node at (1.333333,0.666667) {8}; 
            \fill[black!50] (0.000000,1.500000)--(1.000000,1.500000)--(3.000000,0.500000)--cycle;
            \node at (1.333333,1.166667) {9}; 
            \fill[black!55] (0.000000,1.000000)--(2.000000,0.000000)--(3.000000,0.000000)--cycle;
            \node at (1.666667,0.333333) {10}; 
            \fill[black!60] (0.000000,1.500000)--(2.000000,0.500000)--(3.000000,0.500000)--cycle;
            \node at (1.666667,0.833333) {11};
            \draw (0.0,0.0)--(3.500000,0.0)--(0.0,3.500000)--cycle; 
        \end{tikzpicture}
        \\$k=11$
    \end{minipage}
    \\[1cm]
    \begin{minipage}{0.45\textwidth}
        \centering
        %\begin{tikzpicture}[scale=1.07142857] 
        %\begin{tikzpicture}[scale=1.28571428] 
        %\begin{tikzpicture}[scale=1.39285714] 
        \begin{tikzpicture}[scale=1.33928571] 
            \draw [help lines] (0.0,0.0) grid (3.733332,3.733332); 
            \fill[black!10] (1.399998,0.000000)--(1.399998,1.000000)--(2.399998,0.000000)--cycle;
            \node at (1.733332,0.333333) {1}; 
            \fill[black!15] (0.000000,1.600000)--(0.000000,2.600000)--(1.000000,1.600000)--cycle;
            \node at (0.333333,1.933333) {2}; 
            \fill[black!20] (2.733332,0.000000)--(2.733332,1.000000)--(3.733332,0.000000)--cycle;
            \node at (3.066666,0.333333) {3}; 
            \fill[black!25] (0.399998,0.000000)--(0.399998,1.000000)--(1.399998,1.000000)--cycle;
            \node at (0.733332,0.666667) {4}; 
            \fill[black!30] (0.000000,0.600002)--(0.000000,1.600002)--(1.000000,1.600002)--cycle;
            \node at (0.333333,1.266668) {5}; 
            \fill[black!35] (0.399998,0.000000)--(1.399998,0.000000)--(1.399998,1.000000)--cycle;
            \node at (1.066665,0.333333) {6}; 
            \fill[black!40] (0.399999,1.000000)--(1.399999,1.000000)--(1.399999,2.000000)--cycle;
            \node at (1.066666,1.333333) {7}; 
            \fill[black!45] (1.399998,1.000000)--(1.399998,2.000000)--(2.399998,0.000000)--cycle;
            \node at (1.733332,1.000000) {8}; 
            \fill[black!50] (0.000000,2.600002)--(0.000000,3.600002)--(1.000000,1.600002)--cycle;
            \node at (0.333333,2.600002) {9}; 
            \fill[black!55] (0.000000,3.600000)--(1.000000,1.600000)--(1.000000,2.600000)--cycle;
            \node at (0.666667,2.600000) {10}; 
            \fill[black!60] (1.733330,2.000002)--(2.733330,0.000002)--(2.733330,1.000002)--cycle;
            \node at (2.399997,1.000002) {11}; 
            \fill[black!65] (0.733332,3.000000)--(1.733332,2.000000)--(2.733332,0.000000)--cycle;
            \node at (1.733332,1.666667) {12}; 
            \draw (0.0,0.0)--(3.733332,0.0)--(0.0,3.733332)--cycle; 
        \end{tikzpicture}
        \\$k=12$
    \end{minipage}
    \begin{minipage}{0.45\textwidth}
    \centering
        %\begin{tikzpicture}[scale=1.04347826]
        %\begin{tikzpicture}[scale=1.25217391]
        %\begin{tikzpicture}[scale=1.35652174]
        \begin{tikzpicture}[scale=1.30434782]
            \draw [help lines] (0.0,0.0) grid (3.833333,3.833333); 
            \fill[black!10] (0.000000,1.833333)--(0.000000,2.833333)--(1.000000,1.833333)--cycle;
            \node at (0.333333,2.166667) {1}; 
            \fill[black!14] (0.000000,2.833333)--(0.000000,3.833333)--(1.000000,2.833333)--cycle;
            \node at (0.333333,3.166667) {2}; 
            \fill[black!18] (0.166667,0.000000)--(0.166667,1.000000)--(1.166667,0.000000)--cycle;
            \node at (0.500000,0.333333) {3}; 
            \fill[black!22] (1.500000,1.333333)--(1.500000,2.333333)--(2.500000,1.333333)--cycle;
            \node at (1.833333,1.666667) {4}; 
            \fill[black!26] (0.000000,0.833333)--(0.000000,1.833333)--(1.000000,1.833333)--cycle;
            \node at (0.333333,1.500000) {5}; 
            \fill[black!30] (0.500000,1.333333)--(1.500000,1.333333)--(1.500000,2.333333)--cycle;
            \node at (1.166667,1.666667) {6}; 
            \fill[black!34] (0.000000,2.833333)--(1.000000,1.833333)--(1.000000,2.833333)--cycle;
            \node at (0.666667,2.500000) {7}; 
            \fill[black!38] (0.166667,1.000000)--(1.166667,0.000000)--(2.166667,0.000000)--cycle;
            \node at (1.166667,0.333333) {8}; 
            \fill[black!42] (1.166667,1.000000)--(2.166667,0.000000)--(3.166667,0.000000)--cycle;
            \node at (2.166667,0.333333) {9}; 
            \fill[black!46] (0.166667,1.000000)--(1.166667,1.000000)--(2.166667,0.000000)--cycle;
            \node at (1.166667,0.666667) {10}; 
            \fill[black!50] (1.500000,1.333333)--(2.500000,1.333333)--(3.500000,0.333333)--cycle;
            \node at (2.500000,1.000000) {11}; 
            \fill[black!54] (0.500000,1.333333)--(1.500000,1.333333)--(3.500000,0.333333)--cycle;
            \node at (1.833333,1.000000) {12}; 
            \fill[black!58] (0.500000,1.333333)--(2.500000,0.333333)--(3.500000,0.333333)--cycle;
            \node at (2.166667,0.666667) {13};
            \draw (0.0,0.0)--(3.833333,0.0)--(0.0,3.833333)--cycle; 
        \end{tikzpicture}
        \\$k=13$
    \end{minipage}
    \caption{Optimal $k$-triangle packings for $k=10,\ldots,13$}
    \label{fig: packings2}
\end{figure}
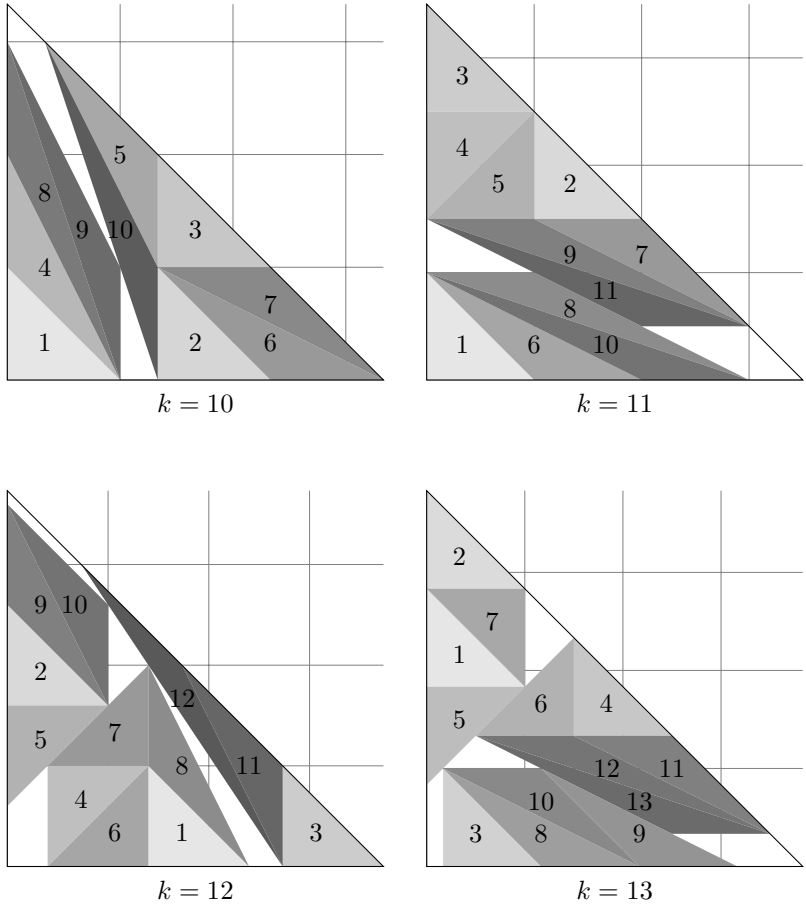

Before we discuss the timing statistics of the algorithm, we take a look at some of the computed packings. Figure~\ref{fig: packings1} and Figure~\ref{fig: packings2} show one exemplary optimal packing for $k=1,\ldots,13$. 
For $k=1,\ldots,12$
we confirm the values of $s_k^{\standardsimplex}$ found by Maley et~al that were presented in Table~\ref{tab: packing-widths}. For  
$k=13$ we can verify their conjecture that the upper bound 
they found is 
indeed optimal.

Our program does not only compute one optimal packing for a given $k$ but is able to detect all multisubsets that allow for an optimal packing. We denote such multisubsets as optimal. In Table~\ref{tab: number of solutions} we report the number of optimal multisubsets for $k = 1,\ldots, 13$.  
For some multisubsets there exist 
more than one optimal packing. 

\begin{table}[ht]
    \centering
    \begin{small}
    \begin{tabular}{|R|R|}
        \hline
        k  &\text{\#Optimal Multisubsets}\\ 
        \hline
        1  &1 \\
        2  &11 \\
        3  &11 \\
        4  &4 \\
        5  &18 \\
        6  &21 \\
        7  &668 \\
        8  &261 \\
        9  &47 \\
        10 &198 \\
        11 &142 \\
        12 &78 \\
        13 &166 \\
        \hline
    \end{tabular}
    \end{small}
    \caption{Number of multisubsets that allow for an optimal $k$-triangle packing for $k = 1,\ldots,13$}
    \label{tab: number of solutions}
\end{table}

An inspection of Table~\ref{tab: number of solutions} raises the natural question why the number of optimal multisubsets increases or decreases significantly for certain $(k,k+1)$ pairs. In order to shine some light on this phenomenon, let us consider sequences for which the optimal container size $s_k^{\standardsimplex}$ does not change, e.g., $k=2,3,4$ or $k=7,8,9$. Assume $k$ and $k-1$ are in the considered sequence.

Henceforth, we abbreviate ``multiset of cardinality $k$'' to ``$k$-multiset''.
From every optimal $k$-multiset with $s$ distinct shapes we obtain $s$ optimal $(k-1)$-multisets (removing one of the first shape or one of the second shape or \dots\ or one of the $s$-th shape).

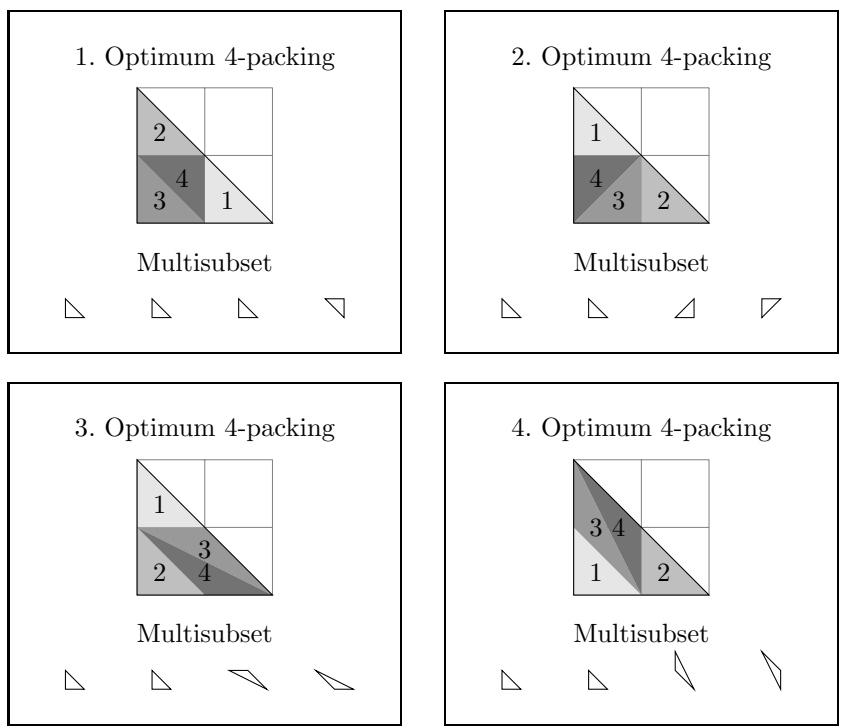
\begin{figure}[!ht]
\begin{center}
\framebox{\parbox{50mm}{
\begin{center}
1.~Optimum 4-packing\\[2mm]
\begin{tikzpicture}[scale=0.9]
\draw [help lines] (0.0,0.0) grid (2.000000,2.000000);
\fill[black!10] (1.000000,0.000000)--(2.000000,0.000000)--(1.000000,1.000000)--cycle;
\node at (1.333333,0.333333) {1};
\fill[black!25] (0.000000,1.000000)--(1.000000,1.000000)--(0.000000,2.000000)--cycle;
\node at (0.333333,1.333333) {2};
\fill[black!40] (0.000000,0.000000)--(1.000000,0.000000)--(0.000000,1.000000)--cycle;
\node at (0.333333,0.333333) {3};
\fill[black!55] (1.000000,0.000000)--(1.000000,1.000000)--(0.000000,1.000000)--cycle;
\node at (0.666667,0.666667) {4};
\draw (0.0,0.0)--(2.000000,0.0)--(0.0,2.000000)--cycle;
\end{tikzpicture}\\[2mm]
Multisubset\\[-6mm]
\vrule height12mm depth0mm width0mm
\hbox to 8mm{\hss\triang{1}{0}{0}{1}\hss}\quad
\hbox to 8mm{\hss\triang{1}{0}{0}{1}\hss}\quad
\hbox to 8mm{\hss\triang{1}{0}{0}{1}\hss}\quad
\hbox to 8mm{\hss\triang{0}{1}{-1}{1}\hss}\quad
\end{center}
}}
\hspace{3.3mm}
\framebox{\parbox{50mm}{
\begin{center}
2.~Optimum 4-packing\\[2mm]
\begin{tikzpicture}[scale=0.9]
\draw [help lines] (0.0,0.0) grid (2.000000,2.000000);
\fill[black!10] (0.000000,1.000000)--(1.000000,1.000000)--(0.000000,2.000000)--cycle;
\node at (0.333333,1.333333) {1};
\fill[black!25] (1.000000,0.000000)--(2.000000,0.000000)--(1.000000,1.000000)--cycle;
\node at (1.333333,0.333333) {2};
\fill[black!40] (0.000000,0.000000)--(1.000000,0.000000)--(1.000000,1.000000)--cycle;
\node at (0.666667,0.333333) {3};
\fill[black!55] (0.000000,0.000000)--(1.000000,1.000000)--(0.000000,1.000000)--cycle;
\node at (0.333333,0.666667) {4};
\draw (0.0,0.0)--(2.000000,0.0)--(0.0,2.000000)--cycle;
\end{tikzpicture}\\[2mm]
Multisubset\\[-6mm]
\vrule height12mm depth0mm width0mm
\hbox to 8mm{\hss\triang{1}{0}{0}{1}\hss}\quad
\hbox to 8mm{\hss\triang{1}{0}{0}{1}\hss}\quad
\hbox to 8mm{\hss\triang{1}{0}{1}{1}\hss}\quad
\hbox to 8mm{\hss\triang{1}{1}{0}{1}\hss}\quad
\end{center}
}}
\end{center}

\begin{center}
\framebox{\parbox{50mm}{
\begin{center}
3.~Optimum 4-packing\\[2mm]
\begin{tikzpicture}[scale=0.9]
\draw [help lines] (0.0,0.0) grid (2.000000,2.000000);
\fill[black!10] (0.000000,1.000000)--(1.000000,1.000000)--(0.000000,2.000000)--cycle;
\node at (0.333333,1.333333) {1};
\fill[black!25] (0.000000,0.000000)--(1.000000,0.000000)--(0.000000,1.000000)--cycle;
\node at (0.333333,0.333333) {2};
\fill[black!40] (2.000000,0.000000)--(1.000000,1.000000)--(0.000000,1.000000)--cycle;
\node at (1.000000,0.666667) {3};
\fill[black!55] (1.000000,0.000000)--(2.000000,0.000000)--(0.000000,1.000000)--cycle;
\node at (1.000000,0.333333) {4};
\draw (0.0,0.0)--(2.000000,0.0)--(0.0,2.000000)--cycle;
\end{tikzpicture}\\[2mm]
Multisubset\\[-6mm]
\vrule height12mm depth0mm width0mm
\hbox to 8mm{\hss\triang{1}{0}{0}{1}\hss}\quad
\hbox to 8mm{\hss\triang{1}{0}{0}{1}\hss}\quad
\hbox to 8mm{\hss\triang{-1}{1}{-2}{1}\hss}\quad
\hbox to 8mm{\hss\triang{1}{0}{-1}{1}\hss}\quad
\end{center}
}}
\hspace{3.3mm}
\framebox{\parbox{50mm}{
\begin{center}
4.~Optimum 4-packing\\[2mm]
\begin{tikzpicture}[scale=0.9]
\draw [help lines] (0.0,0.0) grid (2.000000,2.000000);
\fill[black!10] (0.000000,0.000000)--(1.000000,0.000000)--(0.000000,1.000000)--cycle;
\node at (0.333333,0.333333) {1};
\fill[black!25] (1.000000,0.000000)--(2.000000,0.000000)--(1.000000,1.000000)--cycle;
\node at (1.333333,0.333333) {2};
\fill[black!40] (1.000000,0.000000)--(0.000000,2.000000)--(0.000000,1.000000)--cycle;
\node at (0.333333,1.000000) {3};
\fill[black!55] (1.000000,0.000000)--(1.000000,1.000000)--(0.000000,2.000000)--cycle;
\node at (0.666667,1.000000) {4};
\draw (0.0,0.0)--(2.000000,0.0)--(0.0,2.000000)--cycle;
\end{tikzpicture}\\[2mm]
Multisubset\\[-6mm]
\vrule height12mm depth0mm width0mm
\hbox to 8mm{\hss\triang{1}{0}{0}{1}\hss}\quad
\hbox to 8mm{\hss\triang{1}{0}{0}{1}\hss}\quad
\hbox to 8mm{\hss\triang{-1}{2}{-1}{1}\hss}\quad
\hbox to 8mm{\hss\triang{0}{1}{-1}{2}\hss}\quad
\end{center}
}}
\end{center}
\caption{All optimal multisubsets for $s_4^{\protect\standardsimplex}=2$}
\label{fig: all4packings}
\end{figure}

For $k=4$ we have the 4 optimal multisets shown in Figure~\ref{fig: all4packings}. If we apply the procedure above, the first gives rise to 2 optimal 3-multisets, the other three to 3 optimal 3-multisets, in total exactly 11 multisets, the same that we had calculated before for $k=3$.

Now let us consider the pair $k=2,3$. If we apply the procedure to the 11 optimal multisets for $s_3^{\protect\standardsimplex}=2$, we obtain a total of 24 multisets including duplicates: One of them appears 8 times, six appear twice, four appear once. After removal of the duplicates we obtain the 11 optimal multisets for $k=2$.

While the above examination of the $k=2,3,4$ sequence has been experimental mathematics with pencil and paper, we needed computer assistance for an analogous examination of the $k=7,8,9$ sequence.

\begin{figure}[!ht]
\begin{center}
\framebox{\parbox{80mm}{
\begin{center}
80.~Optimum 7-packing\\[5mm]
\begin{tikzpicture}[scale=1.6]
\draw [help lines] (0.0,0.0) grid (3.000000,3.000000);
\fill[black!10] (0.000000,0.000000)--(1.000000,0.000000)--(0.000000,1.000000)--cycle;
\node at (0.333333,0.333333) {\large 1};
\fill[black!18] (0.666667,0.666667)--(1.666667,0.666667)--(0.666667,1.666667)--cycle;
\node at (1.000000,1.000000) {\large 2};
\fill[black!26] (0.000000,2.000000)--(1.000000,2.000000)--(0.000000,3.000000)--cycle;
\node at (0.333333,2.333333) {\large 3};
\fill[black!34] (2.000000,0.000000)--(3.000000,0.000000)--(2.000000,1.000000)--cycle;
\node at (2.333333,0.333333) {\large 4};
\fill[black!42] (1.000000,0.000000)--(2.000000,0.000000)--(2.000000,1.000000)--cycle;
\node at (1.666667,0.333333) {\large 5};
\fill[black!50] (2.000000,1.000000)--(1.000000,2.000000)--(0.000000,2.000000)--cycle;
\node at (1.000000,1.666667) {\large 6};
\fill[black!58] (1.000000,0.000000)--(0.000000,2.000000)--(0.000000,1.000000)--cycle;
\node at (0.333333,1.000000) {\large 7};
\draw (0.0,0.0)--(3.000000,0.0)--(0.0,3.000000)--cycle;
\end{tikzpicture}\\[3mm]Multisubset\\[-6mm]
\vrule height12mm depth0mm width0mm
\hbox to 6mm{\hss\triang{1}{0}{0}{1}\hss}\quad
\hbox to 6mm{\hss\triang{1}{0}{0}{1}\hss}\quad
\hbox to 6mm{\hss\triang{1}{0}{0}{1}\hss}\quad
\hbox to 6mm{\hss\triang{1}{0}{0}{1}\hss}\quad
\hbox to 6mm{\hss\triang{1}{0}{1}{1}\hss}\quad
\hbox to 6mm{\hss\triang{-1}{1}{-2}{1}\hss}\quad
\hbox to 6mm{\hss\triang{-1}{2}{-1}{1}\hss}\quad
\end{center}
}}
\end{center}
\caption{An example of a nonextendible optimal multisubset for $s_7^{\protect\standardsimplex}=3$}
\label{fig: logo}
\end{figure}
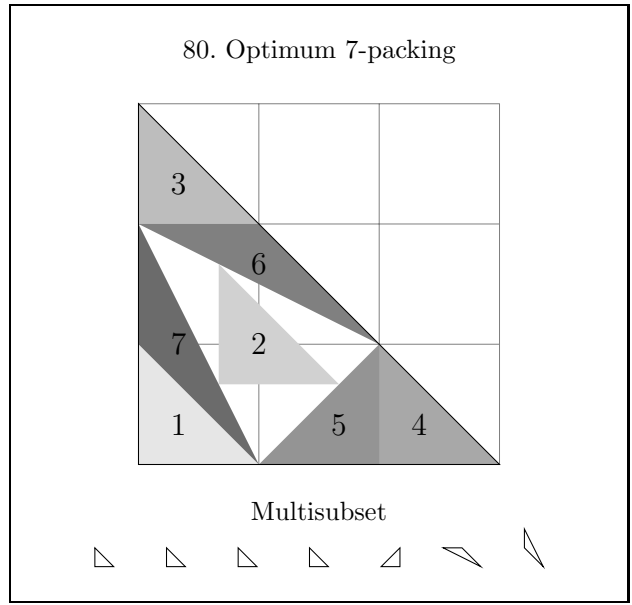

For $k=9$ we have 47 optimal multisets. Application of the procedure gives exactly the 261 multisets that we had calculated before for $k=8$, again with no duplicates.
If we apply the procedure to these 261 8-multisets, we obtain 1361 7-multisets including duplicates. After duplicate removal 646 remain. Beyond these 646, the list of the 668 optimal 7-multisets contains 22 multisets that have no extensions to optimal 8-multisets. 
Since we have a picture of an optimal packing for each optimal multiset we have computed in this study, we could inspect these 22 multisets, one of which is shown in Figure~\ref{fig: logo} along with an optimal packing.\footnote{At this point we can't refrain from telling that a colored version of the picture in Figure~\ref{fig: logo} is the logo of the SFB/TRR 191 `Symplectic Structures in Geometry, Algebra and Dynamics' of the German Science Foundation that has supported our research. When it was chosen in 2016, we had no idea of its rare feature, we just found it visually pleasing.} 

All 22 multisets share the distinct feature that the free spaces in the optimal packing consist of 3 triangles rather than 2 for the other 646 packings whose extendability to 8-packings is visually obvious.

While it is tempting to conjecture that the procedure applied to a full packing (here $k=4$ and $k=9$) gives us the full list of optimum multisets for the next lower $k$ (here 3 and 8), we shall see in Section~\ref{sec: MILP2} that this is not true in general.

Table~\ref{tab: results} shows the timing statistics of our algorithm for $k=1,\ldots,13$. The column labels are:
\begin{itemize}
    \item $k$: the number of triangles to be packed,
    \item $\of{\vector{\abs{\mathcal{S}_k^{\standardsimplex}}}{k}}$: the number of $k$-cardinality multisubsets of the shapelist,
    \item \#I-Calls: the number of calls to the inner optimization procedure,
    \item Avg I-Time: the average cpu time spent in an inner optimization procedure,
    \item Max I-Time: the maximum cpu time spent in an inner optimization procedure,
    \item Total Time: the total cpu time.
\end{itemize}

The number of calls to the inner optimization procedure is significantly smaller than the number of $k$-cardinality multisubsets of the shapelist that would have been considered using complete enumeration. Due to the data base we built up from previous runs for smaller $k$, many calls of the inner optimization procedure have been replaced by simple data base queries. 

\begin{table}[ht]
    \centering
    \begin{small}
    \begin{tabular}{|R|R|R|R|R|R|}
        \hline
        \rule{0pt}{15pt} k
        & \multicolumn{1}{C}{\of{\vector{\abs{\mathcal{S}_k^{\standardsimplex}}}{k}}}
        & \multicolumn{1}{|c}{\#I-Calls} 
        & \multicolumn{1}{|c}{Avg I-Time} 
        & \multicolumn{1}{|c}{Max I-Time} 
        & \multicolumn{1}{|c|}{Total Time}\\[7pt]
        \hline
            1  &1               &1          &0:00:00.00     &0:00:00.00     &0:00:00.00\\
            2  &36              &43         &0:00:00.00     &0:00:00.01     &0:00:00.02\\
            3  &120             &11         &0:00:00.00     &0:00:00.02     &0:00:00.03\\
            4  &330             &11         &0:00:00.01     &0:00:00.09     &0:00:00.14\\
            5  &792             &433        &0:00:00.01     &0:00:00.27     &0:00:03.11\\
            6  &1\,716          &185        &0:00:00.04     &0:00:00.33     &0:00:06.72\\
            7  &657\,800        &255\,158   &0:00:00.00     &0:00:13.30     &0:18:56.90\\\
            8  &2\,220\,075     &263        &0:00:00.14     &0:00:04.46     &0:12:29.10\\
            9  &6\,906\,900     &47         &0:00:00.09     &0:00:01.64     &0:11:59.33\\
            10 &20\,030\,010    &34\,029    &0:00:00.52     &0:02:27.23     &4:56:28.02\\
            11 &54\,627\,300    &43\,187    &0:00:07.67     &0:46:18.30     &92:05:28.83\\
            12 &141\,120\,525   &129\,630   &0:00:09.39     &3:26:34.38     &338:19:31.65\\
            13 &347\,373\,600   &196\,735   &0:00:37.88     &46:59:29.37    &2070:08:23.43\\
        \hline
    \end{tabular}
    \end{small}
    \caption{Timing statistics for the $k$-triangle packing given in the format ``hh:mm:ss" for $k=1,\ldots,13$}
    \label{tab: results}
\end{table}

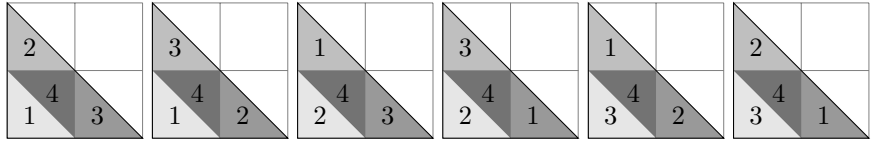
\begin{figure}[ht]
    \centering
    \begin{tikzpicture}[scale = 0.9]
        \draw [help lines] (0.0,0.0) grid (2.000000,2.000000);
        \fill[black!10] (0.000000,0.000000)--(1.000000,0.000000)--(0.000000,1.000000)--cycle;
        \node at (0.333333,0.333333) {1};
        \fill[black!25] (0.000000,1.000000)--(1.000000,1.000000)--(0.000000,2.000000)--cycle;
        \node at (0.333333,1.333333) {2};
        \fill[black!40] (1.000000,0.000000)--(2.000000,0.000000)--(1.000000,1.000000)--cycle;
        \node at (1.333333,0.333333) {3};
        \fill[black!55] (1.000000,0.000000)--(1.000000,1.000000)--(0.000000,1.000000)--cycle;
        \node at (0.666667,0.666667) {4};
        \draw (0.0,0.0)--(2.000000,0.0)--(0.0,2.000000)--cycle;
    \end{tikzpicture}
    \begin{tikzpicture}[scale = 0.9]
        \draw [help lines] (0.0,0.0) grid (2.000000,2.000000);
        \fill[black!10] (0.000000,0.000000)--(1.000000,0.000000)--(0.000000,1.000000)--cycle;
        \node at (0.333333,0.333333) {1};
        \fill[black!25] (0.000000,1.000000)--(1.000000,1.000000)--(0.000000,2.000000)--cycle;
        \node at (0.333333,1.333333) {3};
        \fill[black!40] (1.000000,0.000000)--(2.000000,0.000000)--(1.000000,1.000000)--cycle;
        \node at (1.333333,0.333333) {2};
        \fill[black!55] (1.000000,0.000000)--(1.000000,1.000000)--(0.000000,1.000000)--cycle;
        \node at (0.666667,0.666667) {4};
        \draw (0.0,0.0)--(2.000000,0.0)--(0.0,2.000000)--cycle;
    \end{tikzpicture}
    \begin{tikzpicture}[scale = 0.9]
        \draw [help lines] (0.0,0.0) grid (2.000000,2.000000);
        \fill[black!10] (0.000000,0.000000)--(1.000000,0.000000)--(0.000000,1.000000)--cycle;
        \node at (0.333333,0.333333) {2};
        \fill[black!25] (0.000000,1.000000)--(1.000000,1.000000)--(0.000000,2.000000)--cycle;
        \node at (0.333333,1.333333) {1};
        \fill[black!40] (1.000000,0.000000)--(2.000000,0.000000)--(1.000000,1.000000)--cycle;
        \node at (1.333333,0.333333) {3};
        \fill[black!55] (1.000000,0.000000)--(1.000000,1.000000)--(0.000000,1.000000)--cycle;
        \node at (0.666667,0.666667) {4};
        \draw (0.0,0.0)--(2.000000,0.0)--(0.0,2.000000)--cycle;
    \end{tikzpicture}
    \begin{tikzpicture}[scale = 0.9]
        \draw [help lines] (0.0,0.0) grid (2.000000,2.000000);
        \fill[black!10] (0.000000,0.000000)--(1.000000,0.000000)--(0.000000,1.000000)--cycle;
        \node at (0.333333,0.333333) {2};
        \fill[black!25] (0.000000,1.000000)--(1.000000,1.000000)--(0.000000,2.000000)--cycle;
        \node at (0.333333,1.333333) {3};
        \fill[black!40] (1.000000,0.000000)--(2.000000,0.000000)--(1.000000,1.000000)--cycle;
        \node at (1.333333,0.333333) {1};
        \fill[black!55] (1.000000,0.000000)--(1.000000,1.000000)--(0.000000,1.000000)--cycle;
        \node at (0.666667,0.666667) {4};
        \draw (0.0,0.0)--(2.000000,0.0)--(0.0,2.000000)--cycle;
    \end{tikzpicture}
    \begin{tikzpicture}[scale = 0.9]
        \draw [help lines] (0.0,0.0) grid (2.000000,2.000000);
        \fill[black!10] (0.000000,0.000000)--(1.000000,0.000000)--(0.000000,1.000000)--cycle;
        \node at (0.333333,0.333333) {3};
        \fill[black!25] (0.000000,1.000000)--(1.000000,1.000000)--(0.000000,2.000000)--cycle;
        \node at (0.333333,1.333333) {1};
        \fill[black!40] (1.000000,0.000000)--(2.000000,0.000000)--(1.000000,1.000000)--cycle;
        \node at (1.333333,0.333333) {2};
        \fill[black!55] (1.000000,0.000000)--(1.000000,1.000000)--(0.000000,1.000000)--cycle;
        \node at (0.666667,0.666667) {4};
        \draw (0.0,0.0)--(2.000000,0.0)--(0.0,2.000000)--cycle;
    \end{tikzpicture}
    \begin{tikzpicture}[scale = 0.9]
        \draw [help lines] (0.0,0.0) grid (2.000000,2.000000);
        \fill[black!10] (0.000000,0.000000)--(1.000000,0.000000)--(0.000000,1.000000)--cycle;
        \node at (0.333333,0.333333) {3};
        \fill[black!25] (0.000000,1.000000)--(1.000000,1.000000)--(0.000000,2.000000)--cycle;
        \node at (0.333333,1.333333) {2};
        \fill[black!40] (1.000000,0.000000)--(2.000000,0.000000)--(1.000000,1.000000)--cycle;
        \node at (1.333333,0.333333) {1};
        \fill[black!55] (1.000000,0.000000)--(1.000000,1.000000)--(0.000000,1.000000)--cycle;
        \node at (0.666667,0.666667) {4};
        \draw (0.0,0.0)--(2.000000,0.0)--(0.0,2.000000)--cycle;
    \end{tikzpicture}
    \caption{Symmetric solutions for the multisubset $\offf{T_1, T_1, T_1, T_4}$}
    \label{fig: symmetry}
\end{figure}

One can see that the inner optimization procedure is quite fast on average but a few instances with larger values of $k$ can take very long, especially for $k=13$. We have explored several improvement strategies to the algorithmic approach. The first improvement strategy consists of reducing the occurrence of symmetric solutions for the inner optimization problem. For example, consider the inner optimization instance $\offf{T_1, T_1, T_1, T_4}$. As can be seen in Figure~\ref{fig: symmetry}, we obtain six symmetric optimal solutions for this multisubset.

In general, for each shape with multiplicity $\mu$ we obtain $\mu!$ symmetric solutions. Our approach to reduce the occurrence of symmetric solutions is to add artificial symmetry breaking inequalities to the MILP version of the inner optimization problem. We have examined three different types of symmetry breaking inequalities. All of them try to put an order on the displacement variables $(x_i, y_i)$ belonging to triangles of the same shape. To compare the three symmetry breaking inequality types, we calculated all optimal $k$-triangle packings for $k = 1,\ldots, 12$. Table~\ref{tab: symmetry breaking} shows the computation time of our algorithm combined with the corresponding symmetry breaking constraint. The inequality type numbers are given in the following order: 
\begin{itemize}
    \item Type 0: No symmetry breaking inequality applied.
    \item Type 1: $x_i \leq x_{i+1}$.
    \item Type 2: $x_i+y_i \leq x_{i+1} + y_{i+1}$.
    \item Type 3: $(x_i \leq x_{i+1}) \wedge (x_i = x_{i+1} \Rightarrow y_i \leq y_{i+1})$.
\end{itemize}

\begin{table}[ht]
    \centering
    \begin{small}
   \begin{tabular}{|R|R|R|R|R|}
        \hline
        \multicolumn{1}{|C}{k} & \multicolumn{1}{|C}{\text{Type 0}} & \multicolumn{1}{|C}{\text{Type 1}} & \multicolumn{1}{|C}{\text{Type 2}} & \multicolumn{1}{|C|}{\text{Type 3}}\\
        \hline
        1   &0:00:00.00     &0:00:00.00    &0:00:00.00     &0:00:00.00\\
        2   &0:00:00.04     &0:00:00.04    &0:00:00.02     &0:00:00.03\\
        3   &0:00:00.05     &0:00:00.03    &0:00:00.03     &0:00:00.01\\
        4   &0:00:00.36     &0:00:00.14    &0:00:00.14     &0:00:00.02\\
        5   &0:00:05.67     &0:00:03.16    &0:00:03.11     &0:00:02.88\\
        6   &0:00:12.71     &0:00:09.37    &0:00:06.72     &0:00:08.62\\
        7   &0:28:47.35     &0:17:42.80    &0:18:56.90     &0:14:33.28\\
        8   &0:12:47.89     &0:12:25.75    &0:12:29.10     &0:12:30.34\\
        9   &0:11:57.62     &0:11:57.32    &0:11:59.33     &0:11:59.01\\
        10  &16:04:10.11    &5:33:28.68    &4:56:28.02     &5:29:36.59\\
        11  &1103:05:34.83  &116:21:41.18  &92:05:28.83    &93:30:54.37\\
        12  &4006:41:10.15  &506:53:59.92  &338:19:31.65   &530:56:44.44\\
        \hline
    \end{tabular}
    \end{small}
    \caption{Computation time of the $k$-triangle packing involving symmetry breaking constraints given in the format ``hh:mm:ss" for $k=1,\ldots,12$}
    \label{tab: symmetry breaking}
\end{table}

For $k \geq 10$, one can see that applying symmetry breaking constraints is an enormous improvement to the run-time of the algorithm. For small values of $k$, there is no significant difference between the three symmetry breaking constraints, but for $k = 12$ the clear winner is the symmetry breaking type~2 of the form $x_i+y_i \leq x_{i+1} + y_{i+1}$. For the run producing the computational results presented in Table~\ref{tab: results}, we also chose this symmetry breaking type. 

The second improvement strategy consist in applying a time limit to the inner optimization procedure and changing the GUROBI parameter settings for instances exceeding the given time limit. The GUROBI Optimizer provides a wide variety of parameters that allow  to control the operation of the optimization engines. To find parameter values that improve the performance on our model, we both manually tested different parameter settings and used the built-in automated parameter tuning tool. We found that there are four parameters that speed up the computation time for most instances. Table~\ref{tab: parameter settings} gives an overview of the four parameters together with their default value and the modified value.

\begin{table}[ht]
    \centering
    \strutlongstacks{T}
    \begin{small}
    \begin{tabular}{|l|l|R|R|}
    \hline
        Parameter & \multicolumn{1}{c|}{Description}  & \Centerstack{ Default \\ Value} & \Centerstack{ Modified \\ Value} \\
        \hline
        Presolve    & Controls the presolve level           & -1    & 2\\
        PreDual     & Controls presolve model dualization   & -1    & 1\\ 
        MIPFocus    & MIP solver focus                      & 0     & 2\\
        Heuristics  & Time spent in feasibility heuristics  & 0.05  & 0\\
        \hline
    \end{tabular}
    \end{small}
    \caption{Modified GUROBI parameter setting}
    \label{tab: parameter settings}
\end{table}

For the $k=13$ computation, we reapplied the inner optimization procedure under the modified parameter setting on the instances that exceeded a time limit of 30 seconds. This decreased the total computation time by approximately $38\%$ compared to applying no time limit and using the default parameter setting. The drawback is that we only detect 165 instead of all 166 optimal multisubsets.
In contrast, when applying the modified parameter setting to all inner optimization instances, the computation time increases by approximately $26\%$. The results are summarized in Table~\ref{tab: parameter and time}.

\begin{table}[ht]
    \centering
    \strutlongstacks{T}
    \begin{small}
    \begin{tabular}{|L|l|R|R|R|R|R|R|}
        \hline
        \Centerstack{Time \\Limit} 
        & \Centerstack{Parameter \\Setting}
        & \multicolumn{1}{c|}{Avg I-Time} 
        & \multicolumn{1}{c|}{Total Time}
        & \multicolumn{1}{c|}{\Centerstack{\#Optimal \\Multisubsets}} 
        & \multicolumn{1}{c|}{\Centerstack{\#Open \\Problems}}\\
        \hline
        \infty  &default    &0:00:37.88     &2070:08:23.43  &166    &0\\
        \infty  &modified   &0:00:47.74     &2608:06:35.65  &166    &0\\
        30s     &default    &0:00:18.49     &1001:34:32.45  &160    &1\,650\\
        30s     &modified   &0:00:23.31     &1274:41:01.53  &165    &1\,486\\
        \hline
    \end{tabular}
    \end{small}
    \caption{Timing statistics for the $13$-triangle packing under time limit to the inner optimization procedure combined with modified GUROBI parameter setting given in the format ``hh:mm:ss"}
    \label{tab: parameter and time}
\end{table}

The third improvement strategy is to reformulate the inner optimization problem as a QCQP. Once this formulation is achieved, we can consider the semidefinite relaxation of this program. The hope is that the semidefinite relaxation gives us strong lower bounds on $s$ such that the exact formulation coming from the MILP only needs to be applied if the bounds are not strong enough. This could in turn speed up the branch-and-bound process of the outer optimization problem. We will describe this approach in the next section.

\section{SDP for $k$-Simplex Packing Width of $P^4(1)$}\label{sec: SDP}

\subsection*{Modelling}

We have three different approaches to model the inner optimization problem as a QCQP. We will successively describe each approach and conclude with a comparison of the results.

In the first approach we model the containment condition just like we did in the MILP formulation. Instead of expressing the disjointness condition by a Big-$M$-inequality, we will use the quadratic inequality 
\begin{equation*}
    z^{ij}_f \of{ \alpha^{ij}_f (x_j - x_i) + \beta^{ij}_f (y_j - y_i) - \gamma^{ij}_f } \geq 0.
\end{equation*}
We formulate the binary constraint $z^{ij}_f \in \offf{0,1}$ by the equation $z^{ij}_f(1-z^{ij}_f) = 0$. Since we want at least one of the variables $z^{ij}_f$ to take the value one, we require $z^{ij}_1 + \cdots + z^{ij}_6 \geq 1$. In order to obtain a semidefinite relaxation, we desire the QCQP to be in homogeneous form. To eliminate the linear terms, we introduce a variable $t\in\R$ with $t^2 = 1$ and multiply all linear terms by it. Thus, we obtain the following QCQP formulation of the inner optimization problem.

\begin{Prb}[Inner Optimization Problem - QCQP 1]\label{prb: inner3}
    \begin{align*} 
        \min\;\; & ts \\
        & tx_i \geq K_i & \forall  i \in \offf{1,\ldots,m} \\ 
        & ty_i \geq K'_i & \forall  i \in \offf{1,\ldots,m} \\ 
        & tx_i + ty_i - ts \leq K''_i & \forall  i \in \offf{1,\ldots,m} \\
        & tz^{ij}_1 + \cdots + tz^{ij}_6 \geq 1 & \forall 1 \leq i < j \leq m \\
        & z^{ij}_f \of{ \alpha^{ij}_f (x_j - x_i) + \beta^{ij}_f (y_j - y_i) - t \gamma^{ij}_f } \geq 0 \hspace{-10mm}\\[-1mm]
        && \forall  1 \leq i < j \leq m;; \forall f \in \offf{1,\ldots,6}  \\
        & z^{ij}_f \of{ t - z^{ij}_f } = 0 &  \forall  1 \leq i < j \leq m\;\; \forall f \in \offf{1,\ldots,6} \\
        & t^2 = 1 \\
        & s \in \R, \quad x,y \in \R^m, \quad z \in \R^{6\binom{m}{2}}, \quad t\in\R.\hspace{-10mm}
    \end{align*}
\end{Prb}

The second approach uses the separating hyperplane theorem.
\begin{Thm}[Separating Hyperplane Theorem~\cite{VandenbergheBoyd1996}]\label{thm: separating hyperplane}
    Let $C,D\subseteq\R^n$ be nonempty disjoint convex sets. Then there exist a nonzero vector $a\in\R^n$ and a real number $b$ such that $a^Tx \leq b$ for all $x\in C$ and $a^Tx \geq b$ for all $x\in D$. 
\end{Thm}
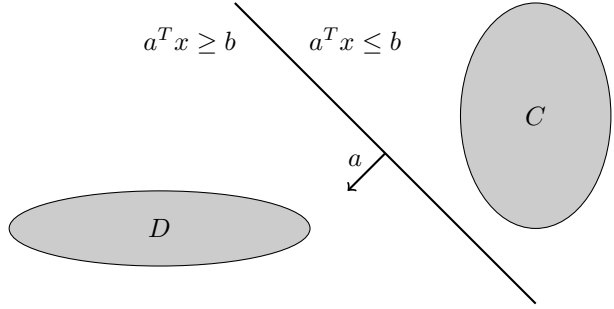
\begin{figure}[ht]
    \centering
    \begin{tikzpicture}
        \node[ellipse, draw, fill=black!20, minimum height = 1cm, minimum width = 4cm] at (-2,0) {$D$};
        \node[ellipse, draw, fill=black!20, minimum height = 3cm, minimum width = 2cm] at (3,1.5) {$C$};
        \draw[thick] (-1,3) -- (3,-1);
        \draw[thick,->] (1,1) -- (0.5,0.5);
        \node at (0.6,0.9) {$a$};
        \node at (0.6,2.5) {$a^Tx \leq b$};
        \node at (-1.6,2.5) {$a^Tx \geq b$};
    \end{tikzpicture}
    \caption{The hyperplane $\{x\in\R^n\mid a^Tx = b\}$ separates the convex sets $C$ and~$D$.}
    \label{fig: separating hyperplane}
\end{figure}
We can directly apply Theorem~\ref{thm: separating hyperplane}, illustrated in Figure~\ref{fig: separating hyperplane}, to the disjointness condition 
\begin{equation*}
    \overline{T_i + t_i} \cap \overline{t_j + T_j} = \emptyset.
\end{equation*}
By Theorem~\ref{thm: separating hyperplane}, the two sets are disjoint if there exist a nonzero vector  $\vector{\alpha^{ij}}{\beta^{ij}}$ and a real number $\gamma^{ij}$ such that  
\begin{align*}
    &\alpha^{ij} x_i                    + \beta^{ij} y_i           &\leq \gamma^{ij}, 
    &\quad\quad\quad \alpha^{ij} x_j         + \beta^{ij} y_j           &\geq \gamma^{ij}, 
    \\
    &\alpha^{ij} (x_i + a_i)            + \beta^{ij} (y_i + b_i)   &\leq \gamma^{ij},
    &\quad\quad\quad \alpha^{ij} (x_j + a_j) + \beta^{ij} (y_j + b_j)   &\geq \gamma^{ij}, 
    \\
    &\alpha^{ij} (x_i + c_i) + \beta^{ij} (y_i + d_i)   &\leq \gamma^{ij}, 
    &\quad\quad\quad \alpha^{ij} (x_j + c_j)            + \beta^{ij} (y_j + d_j)   &\geq \gamma^{ij}.
\end{align*}

Because of convexity, the inequalities hold for the entire set $\overline{t_i+T_i}$ and $\overline{t_j+T_j}$ if they are satisfied by their vertices. To ensure that $\vector{\alpha^{ij}}{\beta^{ij}}$ is nonzero, we impose the additional constraint $\of{\alpha^{ij}}^2 + \of{\beta^{ij}}^2 = 1$. Together with the three containment constraints we obtain the following QCQP. 

\begin{Prb}[Inner Optimization Problem - QCQP 2]\label{prb: inner4}
    \begin{align*}
        \min\;\; & s & & \\
        & tx_i \geq K_i & \forall & i \in \offf{1,\ldots,m} \\ 
        & ty_i \geq K'_i & \forall & i \in \offf{1,\ldots,m} \\ 
        & tx_i + ty_i - ts \leq K''_i & \forall & i \in \offf{1,\ldots,m} \\
        & \of{\alpha^{ij}}^2 + \of{\beta^{ij}}^2 = 1 & \forall & 1 \leq i < j \leq m \\
        & \alpha^{ij} x_i + \beta^{ij} y_i \leq t\gamma^{ij} & \forall & 1 \leq i < j \leq m \\
        & \alpha^{ij} (x_i + ta_i) + \beta^{ij} (y_i + tb_i) \leq t\gamma^{ij} & \forall & 1 \leq i < j \leq m \\
        & \alpha^{ij} (x_i + tc_i) + \beta^{ij} (y_i + td_i) \leq t\gamma^{ij} & \forall & 1 \leq i < j \leq m \\
        & \alpha^{ij} x_j + \beta^{ij} y_j \geq t\gamma^{ij} & \forall & 1 \leq i < j \leq m \\
        & \alpha^{ij} (x_j + ta_j) + \beta^{ij} (y_j + tb_j) \geq t\gamma^{ij} & \forall & 1 \leq i < j \leq m \\
        & \alpha^{ij} (x_j + tc_j) + \beta^{ij} (y_j + td_j) \geq t\gamma^{ij} & \forall & 1 \leq i < j \leq m \\
        & t^2 = 1\\
        & s \in \R, \quad x,y \in \R^m, \quad \alpha,\beta,\gamma \in \R^{\binom{m}{2}}, \quad t\in\R
    \end{align*}
\end{Prb}

The third approach utilizes Farkas Lemma.
\begin{Thm}[Farkas Lemma~\cite{Farkas1902}]
    Let $A\in\R^{m\times n}$ be a matrix and $b\in\R^m$ be a vector. There exists a vector $x\in\R^n$ such that $x\geq 0$ and $Ax=b$ if and only if there does not exist a vector $y\in\R^m$ such that $y^TA\geq 0$ and $y^Tb = -1$.
\end{Thm}

As seen before the disjointness constraint $\overline{t_i + T_i} \cap \overline{t_j + T_j} = \emptyset$ is equivalent to $t_j - t_i \notin \overline{T_i \ominus T_j}$. Previously, we described the Minkowski difference as an intersection of halfspaces. This time we will describe it as the convex hull of its vertices:
\begin{align*}
    T_i \ominus T_j = \operatorname{int} 
    &\left( \operatorname{conv} \left( \left\{ 
    \vector{0}{0}, \vector{a_i}{b_i}, \vector{c_i}{d_i}, \vector{-a_j}{-b_j}, \vector{a_i - a_j}{b_i - b_j}, 
    \right. \right. \right. \\
    &\hspace{1.5cm} \left. \left. \left. \vector{c_i - a_j}{d_i - b_j}, \vector{-c_j}{-d_j}, \vector{a_i - c_j}{b_i - d_j}, \vector{c_i - c_j}{d_i - d_j}
    \right\} \right) \right).
\end{align*}

We know that the difference vector $t_j - t_i$ is not an element of the Minkowski difference $\overline{T_i \ominus T_j}$ if and only if it cannot be written as a convex combination of its vertices. More formally, there do not exist $\lambda^{ij}_1, \ldots, \lambda^{ij}_9 \geq 0$ with $\lambda^{ij}_1 + \cdots + \lambda^{ij}_9 = 1$ such that
\begin{align*}
    \vector{x_j - x_i}{y_j - y_i} = &\lambda^{ij}_1 \vector{0}{0} + \lambda^{ij}_2 \vector{a_i}{b_i} + \lambda^{ij}_3 \vector{c_i}{d_i} + \lambda^{ij}_4 \vector{-a_j}{-b_j} + \lambda^{ij}_5 \vector{a_i - a_j}{b_i - b_j} \\
    &+ \lambda^{ij}_6 \vector{c_i - a_j}{d_i - b_j} + \lambda^{ij}_7 \vector{-c_j}{-d_j} + \lambda^{ij}_8 \vector{a_i - c_j}{b_i - d_j} + \lambda^{ij}_9 \vector{c_i - c_j}{d_i - d_j}.    
\end{align*}

Written in matrix form, this means that there does not exist a vector $\lambda^{ij}\in\R^9$ such that $\lambda^{ij} \geq 0$ 
and
\begin{equation*}
    \begin{bmatrix}
        0 & a_i & c_i & -a_j & a_i-a_j & c_i-a_j & -c_j & a_i-c_j & c_i-c_j \\
        0 & b_i & d_i & -b_j & b_i-b_j & d_i-b_j & -d_j & b_i-d_j & d_i-d_j \\
        1 & 1 & 1 & 1 & 1 & 1 & 1 & 1 & 1
    \end{bmatrix}
    \lambda^{ij} = 
    \begin{bmatrix}
       x_j - x_i \\
       y_j - y_i \\
       1
    \end{bmatrix}\!.
\end{equation*}
Then, by Farkas Lemma, there must exist a vector $\mu^{ij}\in\R^3$ such that 
\begin{equation*}
    \of{\mu^{ij}}^T
    \begin{bmatrix}
        0 & a_i & c_i & -a_j & a_i-a_j & c_i-a_j & -c_j & a_i-c_j & c_i-c_j \\
        0 & b_i & d_i & -b_j & b_i-b_j & d_i-b_j & -d_j & b_i-d_j & d_i-d_j \\
        1 & 1 & 1 & 1 & 1 & 1 & 1 & 1 & 1
    \end{bmatrix}
    \geq 0 
\end{equation*}
and
\begin{equation*}
    \of{\mu^{ij}}^T
    \begin{bmatrix}
       x_j - x_i \\
       y_j - y_i \\
       1
    \end{bmatrix}
    = -1.
\end{equation*}

From the second constraint we derive $\mu^{ij}_3 = - 1 - \mu^{ij}_1 \of{x_j - x_i} - \mu^{ij}_2 \of{y_j - y_i}$. By plugging this expression into the first constraint and renaming $\mu^{ij}_1$ to $\alpha^{ij}$ and $\mu^{ij}_2$ to $\beta^{ij}$, we obtain the following QCQP.

\begin{Prb}[Inner Optimization Problem - QCQP 3]\label{prb: inner5}
    \begin{align*}
        \min\;\; & s & & \\
        & tx_i \geq K_i & \forall & i \in \offf{1,\ldots,m} \\ 
        & ty_i \geq K'_i & \forall & i \in \offf{1,\ldots,m} \\ 
        & tx_i + ty_i - ts \leq K''_i & \forall & i \in \offf{1,\ldots,m} \\
        & - \alpha^{ij} \of{x_j - x_i} - \beta^{ij} \of{y_j - y_i} \geq 1 & \forall & 1 \leq i < j \leq m \\
        & \alpha^{ij} \of{ta_i - x_j + x_i} + \beta^{ij} \of{tb_i - y_j + y_i} \geq 1 & \forall & 1 \leq i < j \leq m \\
        & \alpha^{ij} \of{tc_i - x_j + x_i} + \beta^{ij} \of{td_i - y_j + y_i} \geq 1 & \forall & 1 \leq i < j \leq m \\
        & \alpha^{ij} \of{-ta_j - x_j + x_i} + \beta^{ij} \of{-tb_j - y_j + y_i} \geq 1 & \forall & 1 \leq i < j \leq m \\
        & \alpha^{ij} \of{t(a_i - a_j) - x_j + x_i} + \beta^{ij} \of{t(b_i - b_j) - y_j + y_i} \geq 1 & \forall & 1 \leq i < j \leq m \\
        & \alpha^{ij} \of{t(c_i - a_j) - x_j + x_i} + \beta^{ij} \of{t(d_i - b_j) - y_j + y_i} \geq 1 & \forall & 1 \leq i < j \leq m \\
        & \alpha^{ij} \of{-tc_j - x_j + x_i} + \beta^{ij} \of{-td_j - y_j + y_i} \geq 1 & \forall & 1 \leq i < j \leq m \\
        & \alpha^{ij} \of{t(a_i - c_j) - x_j + x_i} + \beta^{ij} \of{t(b_i - d_j) - y_j + y_i} \geq 1 & \forall & 1 \leq i < j \leq m \\
        & \alpha^{ij} \of{t(c_i - c_j) - x_j + x_i} + \beta^{ij} \of{t(d_i - d_j) - y_j + y_i} \geq 1 & \forall & 1 \leq i < j \leq m \\
        & t^2 = 1\\
        & s \in \R, \quad x,y \in \R^m, \quad \alpha,\beta \in \R^{\binom{m}{2}}, \quad t\in\R
    \end{align*}
\end{Prb}

\subsection*{Implementation and Results}

Table~\ref{tab: sdp value comparison} shows the optimal value of the MILP and the three SDPs and Table~\ref{tab: sdp run time comparison} shows the corresponding computation times of the four optimization problems. The instances shown in the first column where randomly chosen. For solving the SDPs, we call the MOSEK Optimization Software Version 8.1.

\begin{table}[ht]
    \centering
    \strutlongstacks{T}
    \begin{small}
    \begin{tabular}{|C|C|C|C|C|}
        \hline
        \multirow{2}{*}{\Centerstack{Random \\Instance}} & \multicolumn{4}{c|}{Optimal Value}\\
        \cline{2-5}
        & \text{MILP} & \text{SDP 1} & \text{SDP 2} & \text{SDP 3} \\
        \hline
        1&3.89 &3.00 &1.00 &1.00 \\
        2&3.92 &3.00 &2.00 &2.00 \\
        3&3.94 &2.00 &2.00 &2.00 \\
        4&3.83 &2.00 &2.00 &2.00 \\
        5&3.81 &2.00 &2.00 &2.00 \\
        6&3.80 &2.00 &2.00 &2.00 \\
        7&3.80 &2.00 &2.00 &2.00 \\
        8&3.83 &2.00 &2.00 &2.00 \\
        9&4.00 &3.00 &3.00 &3.00 \\
        10&3.92 &3.00 &2.00 &2.00 \\
        \hline
    \end{tabular}
    \end{small}
    \caption{Optimal values of the MILP and the three SDPs for different inner optimization instances}
    \label{tab: sdp value comparison}
\end{table}

\begin{table}[ht]
    \centering
    \begin{small}
    \begin{tabular}{|C|C|C|C|}
        \hline
        \multicolumn{4}{|c|}{Computation Time}\\
        \hline
        \text{MILP} & \text{SDP 1} & \text{SDP 2} & \text{SDP 3}\\
        \hline
        0:13:26.6 &0:00:03.7 &0:00:00.4 &0:00:00.3 \\
        0:07:06.2 &0:00:04.2 &0:00:00.3 &0:00:00.3 \\
        1:13:16.2 &0:00:08.0 &0:00:00.7 &0:00:00.4 \\
        1:08:18.9 &0:00:07.1 &0:00:00.7 &0:00:00.5 \\
        0:56:57.6 &0:00:06.1 &0:00:00.6 &0:00:00.6 \\
        0:46:36.9 &0:00:06.5 &0:00:00.6 &0:00:00.4 \\
        0:45:26.9 &0:00:06.8 &0:00:00.7 &0:00:00.5 \\
        0:09:21.5 &0:00:07.5 &0:00:00.9 &0:00:00.4 \\
        1:17:29.6 &0:00:12.4 &0:00:00.9 &0:00:00.6 \\
        0:19:24.2 &0:00:17.3 &0:00:01.2 &0:00:00.8 \\
        \hline
    \end{tabular}
    \end{small}
    \caption{Computation time of the MILP and the three SDPs for different inner optimization instances given in the format ``hh:mm:ss". The instances are in the same order as in Table~\ref{tab: sdp value comparison}.}
    \label{tab: sdp run time comparison}
\end{table}

For the first two and the last inner optimization instances, the optimal values of the three semidefinite relaxations differ. For the remaining inner optimization instances the optimal values coincide. The computation time of the first semidefinite relaxation is slightly larger than the computation time of the other two relaxations but much faster than the computation time of the MILP. However, the gap between the exact value of the inner optimization problem obtained by the MILP and the best lower bound obtained by the semidefinite relaxations is too large for a successful application in our branch and bound framework.

\section{MILP for $k$-Simplex Packing Width of $P^6(1)$}\label{sec: MILP2}

In this section we will extend the MILP approach to the next higher dimension and compute the $k$-simplex packing width of the six-dimensional open prism $P^6(1) = \standardsimplex^3(1) \times \T^3$. The outer and inner optimization problem are defined analogously to the two-dimensional setting. The number of $k$-cardinality multisubsets of the shapelist $\mathcal{S}_k^{\standardtetrahedron}$ for $k = 1,\ldots,8$ is shown in Table~\ref{tab: 6multisubsets}.

\begin{table}[ht]
    \centering
    \begin{small}
    \begin{tabular}{|R|R|R|}                                 
        \hline
        \rule{0pt}{15pt} k & \multicolumn{1}{C}{\abs{\mathcal{S}_k^{\standardtetrahedron}}} & \multicolumn{1}{|C|}{\of{\vector{\abs{\mathcal{S}_k^{\standardtetrahedron}}}{k}}} \\[7pt]
        \hline
        1 & 1 & 1\\
        2 & 73 & 2\,701\\
        3 & 73 & 67\,525 \\
        4 & 73 & 1\,282\,975\\
        5 & 73 & 19\,757\,815\\
        6 & 73 & 256\,851\,595\\
        7 & 73 & 2\,898\,753\,715\\
        8 & 73 & 28\,987\,537\,150\\
        \hline
    \end{tabular}
    \end{small}
    \caption{Number of $k$-cardinality multisubsets of the shapelists for $k=1,\ldots,8$}
    \label{tab: 6multisubsets}
\end{table}

Due to the cardinality of the shapelist, the number of multisubsets is quite large even for small numbers of $k$. However, we were able to compute $k$-tetrahedron packings for $k = 1,\ldots,8$ by employing the same branch-and-bound strategy as in the two-dimensional setting.

\begin{figure}[ht]
    \centering
        \begin{minipage}{0.42\textwidth}
        \centering
        \begin{tikzpicture}[scale=1]                                                    
            \draw[->] (0,0) -- (3.000000,0,0) node[right] {$x$};                                  
            \draw[->] (0,0) -- (0,3.000000,0) node[above] {$y$};                                  
            \draw[->] (0,0) -- (0,0,3.000000) node[below left] {$z$}; 							   
            \begin{scope}[shift={(0.000000,1.000000,0.000000)}]                                           
            \vertices{0}{0}{1}{0}{1}{0}{1}{0}{0}                               
            \edges{black!10}                                                            
            \end{scope}                                                                 
            \node at (0.333333,1.333333,0.333333) {1}; 
            \begin{scope}[shift={(0.000000,0.000000,0.000000)}]                                           
            \vertices{0}{0}{1}{0}{1}{0}{1}{0}{0}                               
            \edges{black!18}                                                            
            \end{scope}                                                                 
            \node at (0.333333,0.333333,0.333333) {2}; 
            \begin{scope}[shift={(0.000000,0.000000,1.000000)}]                                           
            \vertices{0}{0}{1}{0}{1}{0}{1}{0}{0}                               
            \edges{black!26}                                                            
            \end{scope}                                                                 
            \node at (0.333333,0.333333,1.333333) {3}; 
            \begin{scope}[shift={(1.000000,0.000000,0.000000)}]                                           
            \vertices{0}{0}{1}{0}{1}{0}{1}{0}{0}                               
            \edges{black!34}                                                            
            \end{scope}                                                                 
            \node at (1.333333,0.333333,0.333333) {4}; 
            \begin{scope}[shift={(0.000000,0.000000,1.000000)}]                                           
            \vertices{0}{1}{-1}{0}{1}{0}{1}{0}{-1}                               
            \edges{black!42}                                                            
            \end{scope}                                                                 
            \node at (0.333333,0.666667,0.333333) {5}; 
            \begin{scope}[shift={(0.000000,1.000000,0.000000)}]                                           
            \vertices{0}{0}{1}{1}{-1}{0}{1}{0}{0}                               
            \edges{black!50}                                                            
            \end{scope}                                                                 
            \node at (0.666667,0.666667,0.333333) {6}; 
            \begin{scope}[shift={(0.000000,1.000000,1.000000)}]                                           
            \vertices{1}{-1}{-1}{1}{-1}{0}{1}{0}{-1}                               
            \edges{black!58}                                                            
            \end{scope}                                                                 
            \node at (1.000000,0.333333,0.333333) {7}; 
        \end{tikzpicture}
        \\[-4mm]$k=7$\vspace{4mm}
        \end{minipage}
        \begin{minipage}{0.42\textwidth}
        \centering
        \begin{tikzpicture}[scale=1]                                                    
            \draw[->] (0,0) -- (3.000000,0,0) node[right] {$x$};                                  
            \draw[->] (0,0) -- (0,3.000000,0) node[above] {$y$};                                  
            \draw[->] (0,0) -- (0,0,3.000000) node[below left] {$z$}; 							   
            \begin{scope}[shift={(0.000000,0.000000,1.000000)}]                                           
            \vertices{0}{0}{1}{0}{1}{0}{1}{0}{0}                               
            \edges{black!10}                                                            
            \end{scope}                                                                 
            \node at (0.333333,0.333333,1.333333) {1}; 
            \begin{scope}[shift={(1.000000,0.000000,0.000000)}]                                           
            \vertices{0}{0}{1}{0}{1}{0}{1}{0}{0}                               
            \edges{black!17}                                                            
            \end{scope}                                                                 
            \node at (1.333333,0.333333,0.333333) {2}; 
            \begin{scope}[shift={(0.000000,0.000000,0.000000)}]                                           
            \vertices{0}{0}{1}{0}{1}{0}{1}{0}{1}                               
            \edges{black!24}                                                            
            \end{scope}                                                                 
            \node at (0.333333,0.333333,0.666667) {3}; 
            \begin{scope}[shift={(0.000000,0.000000,1.000000)}]                                           
            \vertices{0}{1}{-1}{0}{2}{-1}{1}{0}{0}                               
            \edges{black!31}                                                            
            \end{scope}                                                                 
            \node at (0.333333,1.000000,0.333333) {4}; 
            \begin{scope}[shift={(0.000000,0.000000,1.000000)}]                                           
            \vertices{0}{1}{0}{0}{2}{-1}{1}{0}{0}                               
            \edges{black!38}                                                            
            \end{scope}                                                                 
            \node at (0.333333,1.000000,0.666667) {5}; 
            \begin{scope}[shift={(0.000000,0.000000,0.000000)}]                                           
            \vertices{0}{1}{0}{1}{0}{0}{1}{0}{1}                               
            \edges{black!45}                                                            
            \end{scope}                                                                 
            \node at (0.666667,0.333333,0.333333) {6}; 
            \begin{scope}[shift={(0.000000,1.000000,0.000000)}]                                           
            \vertices{0}{1}{0}{1}{-1}{1}{1}{0}{0}                               
            \edges{black!52}                                                            
            \end{scope}                                                                 
            \node at (0.666667,1.000000,0.333333) {7}; 
            \begin{scope}[shift={(0.000000,1.000000,0.000000)}]                                           
            \vertices{1}{-1}{0}{1}{-1}{1}{1}{0}{0}                               
            \edges{black!59}                                                            
            \end{scope}                                                                 
            \node at (1.000000,0.333333,0.333333) {8}; 
        \end{tikzpicture}
        \\[-4mm]$k=8$\vspace{4mm}
        \end{minipage}
    \caption{Optimal $k$-tetrahedron packings for $k = 7,8$}
    \label{fig: 6packings2}
\end{figure}
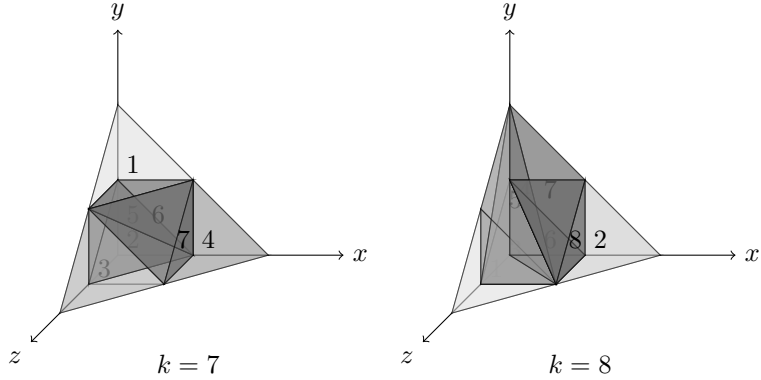

We found that $s_1^{\standardtetrahedron} = 1$ and $s_k^{\standardtetrahedron} = 2$ for $k = 2,\ldots,8$. Figure~\ref{fig: 6packings1} and Figure~\ref{fig: 6packings2} show one exemplary optimal packing for $k=1,\ldots,8$.

\begin{figure}[!ht]
    \centering
    \begin{minipage}{0.42\textwidth}
        \centering
        \begin{tikzpicture}[scale=1.5]                                                    
            \draw[->] (0,0) -- (2.000000,0,0) node[right] {$x$};                                  
            \draw[->] (0,0) -- (0,2.000000,0) node[above] {$y$};                                  
            \draw[->] (0,0) -- (0,0,2.000000) node[below left] {$z$};                             
            \begin{scope}[shift={(0.000000,0.000000,0.000000)}]                                           
            \vertices{0}{0}{1}{0}{1}{0}{1}{0}{0}                               
            \edges{black!10}                                                            
            \end{scope}                                                                 
            \node at (0.333333,0.333333,0.333333) {1}; 
        \end{tikzpicture}
        \\[-4mm]$k=1$\vspace{6mm}
    \end{minipage}
    \begin{minipage}{0.42\textwidth}
        \centering
        \begin{tikzpicture}[scale=1]                                                    
            \draw[->] (0,0) -- (3.000000,0,0) node[right] {$x$};                                  
            \draw[->] (0,0) -- (0,3.000000,0) node[above] {$y$};                                  
            \draw[->] (0,0) -- (0,0,3.000000) node[below left] {$z$}; 							   
            \begin{scope}[shift={(0.000000,0.000000,2.000000)}]                                           
            \vertices{0}{1}{-2}{0}{1}{-1}{1}{0}{-2}                               
            \edges{black!10}                                                            
            \end{scope}                                                                 
            \node at (0.333333,0.666667,0.333333) {1}; 
            \begin{scope}[shift={(0.000000,2.000000,0.000000)}]                                           
            \vertices{1}{-2}{0}{1}{-2}{1}{1}{-1}{0}                               
            \edges{black!40}                                                            
            \end{scope}                                                                 
            \node at (1.000000,0.333333,0.333333) {2}; 
        \end{tikzpicture}
        \\[-4mm]$k=2$\vspace{6mm}
        \end{minipage}
        \\
        \begin{minipage}{0.42\textwidth}
        \centering
        \begin{tikzpicture}[scale=1]                                                    
            \draw[->] (0,0) -- (3.000000,0,0) node[right] {$x$};                                  
            \draw[->] (0,0) -- (0,3.000000,0) node[above] {$y$};                                  
            \draw[->] (0,0) -- (0,0,3.000000) node[below left] {$z$}; 							   
            \begin{scope}[shift={(0.000000,0.000000,1.000000)}]                                           
            \vertices{0}{1}{-1}{0}{1}{0}{1}{0}{0}                               
            \edges{black!10}                                                            
            \end{scope}                                                                 
            \node at (0.333333,0.666667,0.666667) {1}; 
            \begin{scope}[shift={(0.000000,0.000000,1.000000)}]                                           
            \vertices{0}{1}{-1}{1}{0}{-1}{1}{0}{0}                               
            \edges{black!30}                                                            
            \end{scope}                                                                 
            \node at (0.666667,0.333333,0.333333) {2}; 
            \begin{scope}[shift={(0.000000,1.000000,0.000000)}]                                           
            \vertices{0}{1}{0}{1}{-1}{0}{1}{-1}{1}                               
            \edges{black!50}                                                            
            \end{scope}                                                                 
            \node at (0.666667,0.666667,0.333333) {3}; 
        \end{tikzpicture}
        \\[-4mm]$k=3$\vspace{6mm}
        \end{minipage}
        \begin{minipage}{0.42\textwidth}
        \centering
        \begin{tikzpicture}[scale=1]                                                    
            \draw[->] (0,0) -- (3.000000,0,0) node[right] {$x$};                                  
            \draw[->] (0,0) -- (0,3.000000,0) node[above] {$y$};                                  
            \draw[->] (0,0) -- (0,0,3.000000) node[below left] {$z$}; 							   
            \begin{scope}[shift={(0.000000,0.000000,0.000000)}]                                           
            \vertices{0}{0}{1}{0}{1}{1}{1}{0}{0}                               
            \edges{black!10}                                                            
            \end{scope}                                                                 
            \node at (0.333333,0.333333,0.666667) {1}; 
            \begin{scope}[shift={(0.000000,0.000000,1.000000)}]                                           
            \vertices{0}{0}{1}{0}{1}{0}{1}{0}{-1}                               
            \edges{black!25}                                                            
            \end{scope}                                                                 
            \node at (0.333333,0.333333,1.000000) {2}; 
            \begin{scope}[shift={(0.000000,0.000000,2.000000)}]                                           
            \vertices{0}{1}{-1}{1}{0}{-2}{1}{0}{-1}                               
            \edges{black!40}                                                            
            \end{scope}                                                                 
            \node at (0.666667,0.333333,0.666667) {3}; 
            \begin{scope}[shift={(0.000000,1.000000,1.000000)}]                                           
            \vertices{0}{1}{-1}{1}{-1}{-1}{1}{0}{-1}                               
            \edges{black!55}                                                            
            \end{scope}                                                                 
            \node at (0.666667,1.000000,0.000000) {4}; 
        \end{tikzpicture}
        \\[-4mm]$k=4$\vspace{6mm}
        \end{minipage}
        \\
        \begin{minipage}{0.42\textwidth}
        \centering
        \begin{tikzpicture}[scale=1]                                                    
            \draw[->] (0,0) -- (3.000000,0,0) node[right] {$x$};                                  
            \draw[->] (0,0) -- (0,3.000000,0) node[above] {$y$};                                  
            \draw[->] (0,0) -- (0,0,3.000000) node[below left] {$z$}; 							   
            \begin{scope}[shift={(0.000000,0.000000,0.000000)}]                                           
            \vertices{0}{0}{1}{0}{1}{0}{1}{0}{0}                               
            \edges{black!10}                                                            
            \end{scope}                                                                 
            \node at (0.333333,0.333333,0.333333) {1}; 
            \begin{scope}[shift={(0.000000,1.000000,0.000000)}]                                           
            \vertices{0}{0}{1}{0}{1}{0}{1}{0}{0}                               
            \edges{black!22}                                                            
            \end{scope}                                                                 
            \node at (0.333333,1.333333,0.333333) {2}; 
            \begin{scope}[shift={(0.000000,0.000000,1.000000)}]                                           
            \vertices{0}{0}{1}{0}{1}{0}{1}{0}{0}                               
            \edges{black!34}                                                            
            \end{scope}                                                                 
            \node at (0.333333,0.333333,1.333333) {3}; 
            \begin{scope}[shift={(0.000000,0.000000,1.000000)}]                                           
            \vertices{0}{1}{-1}{0}{1}{0}{1}{1}{-1}                               
            \edges{black!46}                                                            
            \end{scope}                                                                 
            \node at (0.333333,1.000000,0.333333) {4}; 
            \begin{scope}[shift={(0.000000,0.000000,1.000000)}]                                           
            \vertices{0}{1}{-1}{1}{0}{-1}{2}{0}{-1}                               
            \edges{black!58}                                                            
            \end{scope}                                                                 
            \node at (1.000000,0.333333,0.000000) {5}; 
        \end{tikzpicture}
        \\[-4mm]$k=5$\vspace{6mm}
        \end{minipage}
        \begin{minipage}{0.42\textwidth}
        \centering
        \begin{tikzpicture}[scale=1]                                                    
            \draw[->] (0,0) -- (3.000000,0,0) node[right] {$x$};                                  
            \draw[->] (0,0) -- (0,3.000000,0) node[above] {$y$};                                  
            \draw[->] (0,0) -- (0,0,3.000000) node[below left] {$z$}; 							   
            \begin{scope}[shift={(0.000000,1.000000,0.000000)}]                                           
            \vertices{0}{0}{1}{0}{1}{0}{1}{0}{0}                               
            \edges{black!10}                                                            
            \end{scope}                                                                 
            \node at (0.333333,1.333333,0.333333) {1}; 
            \begin{scope}[shift={(0.000000,0.000000,1.000000)}]                                           
            \vertices{0}{0}{1}{0}{1}{0}{1}{0}{0}                               
            \edges{black!20}                                                            
            \end{scope}                                                                 
            \node at (0.333333,0.333333,1.333333) {2}; 
            \begin{scope}[shift={(1.000000,0.000000,0.000000)}]                                           
            \vertices{0}{0}{1}{0}{1}{0}{1}{0}{0}                               
            \edges{black!30}                                                            
            \end{scope}                                                                 
            \node at (1.333333,0.333333,0.333333) {3}; 
            \begin{scope}[shift={(0.000000,0.000000,0.000000)}]                                           
            \vertices{0}{1}{0}{0}{1}{1}{1}{0}{0}                               
            \edges{black!40}                                                            
            \end{scope}                                                                 
            \node at (0.333333,0.666667,0.333333) {4}; 
            \begin{scope}[shift={(0.000000,0.000000,1.000000)}]                                           
            \vertices{0}{1}{0}{1}{0}{-1}{1}{0}{0}                               
            \edges{black!50}                                                            
            \end{scope}                                                                 
            \node at (0.666667,0.333333,0.666667) {5}; 
            \begin{scope}[shift={(0.000000,1.000000,0.000000)}]                                           
            \vertices{0}{0}{1}{1}{-1}{0}{1}{0}{0}                               
            \edges{black!60}                                                            
            \end{scope}                                                                 
            \node at (0.666667,0.666667,0.333333) {6}; 
        \end{tikzpicture}
        \\[-4mm]$k=6$\vspace{6mm}
        \end{minipage}
    \caption{Optimal $k$-tetrahedron packings for $k = 1,\ldots,6$}
    \label{fig: 6packings1}
\end{figure}
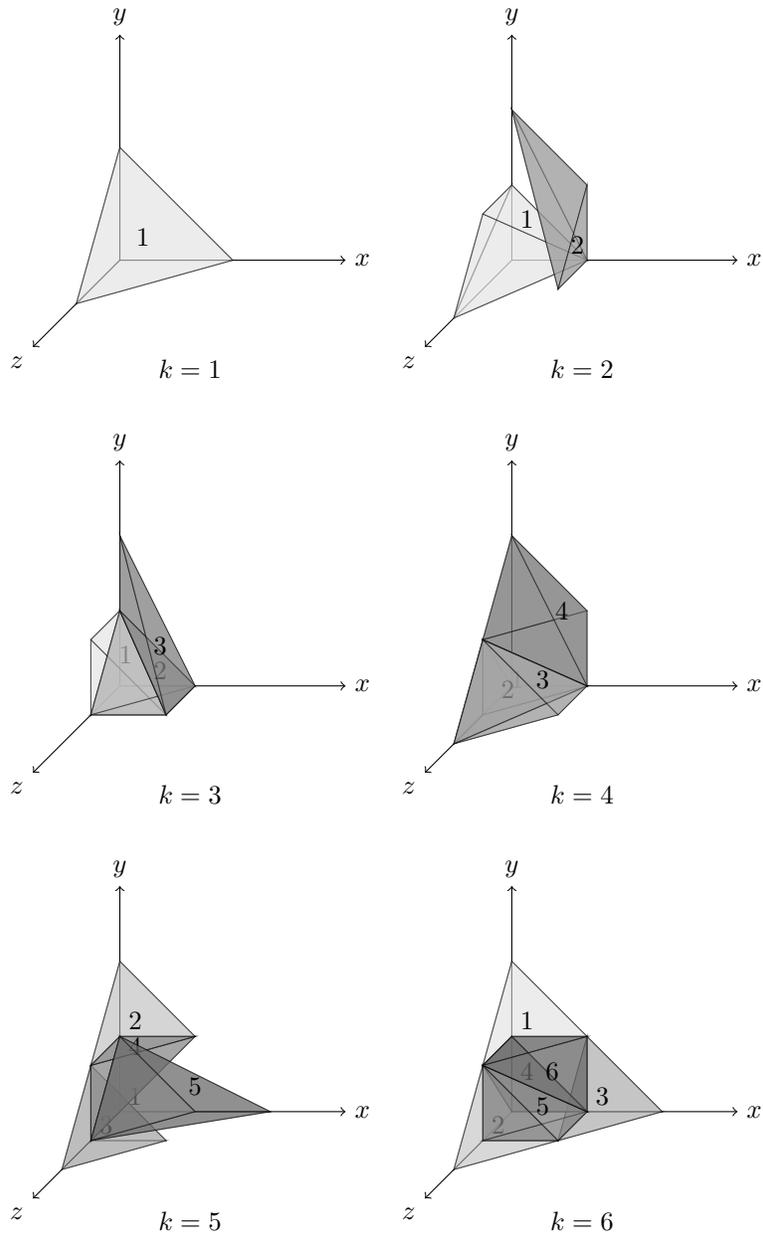

As in the two-dimensional setting,  our program detects all multisubsets that allow for an optimal packing for a given $k$. Table~\ref{tab: 6number of solutions} shows the number of optimal multisubsets for $k = 1,\ldots, 8$.  

\begin{table}[!ht]
    \centering
    \begin{small}
    \begin{tabular}{|R|R|}
        \hline
        k  &\text{\#Optimal Multisubsets}\\ 
        \hline
        1  &1 \\
        2  &1\,123 \\
        3  &4\,871 \\
        4  &9\,914 \\
        5  &11\,709 \\
        6  &8\,075 \\
        7  &2\,899 \\
        8  &408 \\
        \hline
    \end{tabular}
    \end{small}
    \caption{Number of multisubsets that allow for an optimal $k$-tetrahedron packing for $k = 1,\ldots,8$}
    \label{tab: 6number of solutions}
\end{table}

There are far more optimal multisubsets than in the two-dimensional setting due to the greater number of feasible multisubsets. 

In order to contribute to understanding the sequence of numbers in Table~\ref{tab: 6number of solutions}, we have applied the derivation procedure of Section~\ref{sec: MILP1} for $k=8,7,\ldots,2$, all with the same $s_k^{\standardtetrahedron}=2$. 

In Table~\ref{tab: statsderived} we present for each pair $(k,k-1)$ the total number of derived multisets and the number of remaining multisets after duplicate removal, along with the computed number of optimal $(k-1)$-multisets. The last column gives the number of $(k-1)$-multisets that are not extendable to $k$-multisets. As stated in Section~\ref{sec: MILP1}, it is not true that for a full $k$-packing all $(k-1)$-multisets are obtained by the derivation procedure: 64 of the 2899 optimal 7-multisets have no extension to an optimal 8-multiset.

\begin{table}[ht]
\begin{center}
\begin{tabular}{|c|r|r|r|r|}
\hline
&\multicolumn{2}{c|}{number of}&\multicolumn{2}{c|}{}\\
&\multicolumn{2}{c|}{derived multisubsets}&\multicolumn{2}{c|}{number of}\\
\cline{2-5}
&\multicolumn{1}{c|}{including}&\multicolumn{1}{c|}{without}&\multicolumn{1}{c|}{optimal}&\multicolumn{1}{c|}{nonextendable}\\
$k\rightarrow(k-1)$&\multicolumn{1}{c|}{duplicates}&\multicolumn{1}{c|}{duplicates}&\multicolumn{1}{c|}{multisubsets}&\multicolumn{1}{c|}{multisubsets}\\
\hline
$8\rightarrow7$&  2835&  2835&  2899& 64\\
$7\rightarrow6$& 17755&  8071&  8075&  4\\
$6\rightarrow5$& 43353& 11709& 11709&  0\\
$5\rightarrow4$& 54429&  9914&  9914&  0\\
$4\rightarrow3$& 38459&  4871&  4871&  0\\
$3\rightarrow2$& 14539&  1123&  1123&  0\\
\hline
\end{tabular}
\end{center}
    \caption{Statistics on all derived vs.\ all optimal multisets}
    \label{tab: statsderived}
\end{table}

We hope that our experimental findings can serve as a basis for a rigorous mathematical analysis.

Table~\ref{tab: 6results} shows the timing statistics of our algorithm for $k=1,\ldots,8$. The column labels are the same as in the two-dimensional case.

\begin{table}[ht]
    \centering
    \begin{small}
    \begin{tabular}{|R|R|R|R|R|R|}
        \hline
        \rule{0pt}{15pt} k
        & \multicolumn{1}{C}{\of{\vector{\abs{\mathcal{S}_k^{\standardtetrahedron}}}{k}}}
        & \multicolumn{1}{|C}{\#\text{I-Calls}} 
        & \multicolumn{1}{|C}{\text{Avg I-Time}} 
        & \multicolumn{1}{|C}{\text{Max I-Time}} 
        & \multicolumn{1}{|C|}{\text{Total Time}}\\[7pt]
        \hline
            1  &1                   &1              &0:00:00.00     &0:00:00.00     &0:00:00.00\\
            2  &2\,701              &2\,773         &0:00:00.00     &0:00:00.07     &0:00:01.90\\
            3  &67\,525             &4\,871         &0:00:00.00     &0:00:00.04     &0:00:06.59\\
            4  &1\,282\,975         &9\,914         &0:00:00.00     &0:00:00.14     &0:00:25.83\\
            5  &19\,757\,815        &13\,118        &0:00:00.00     &0:00:10.15     &0:06:04.31\\
            6  &256\,851\,595       &8\,075         &0:00:00.01     &0:00:00.09     &0:01:01.06\\
            7  &2\,898\,753\,715    &2\,899         &0:00:00.01     &0:00:00.35     &0:00:39.66\\
            8  &28\,987\,537\,150   &408            &0:00:00.01     &0:00:00.08     &0:00:12.77\\
        \hline
    \end{tabular}
    \end{small}
    \caption{Timing statistics for the $k$-tetrahedron packing given in the format ``hh:mm:ss" for $k=1,\ldots,8$}
    \label{tab: 6results}
\end{table}

The inner optimization procedure is very fast on all instances. For computing $k$-tetrahedron packings for $k \geq 9$, the difficulty rather consists in the high number of feasible multisubsets. For $k = 9$ the cardinality of the shapelist increases to $854$ and the number of multisubsets thereof increases to \[694\,392\,240\,786\,929\,755\,070 \approx 7 \times 10^{20}.\] Therefore, computing $s_9^{\standardtetrahedron}$ seems out of reach. Nevertheless, instead of working with the shapelist $S^{\standardtetrahedron}_9$, one might consider the smaller shapelist $S^{\standardtetrahedron}_8$ to compute upper bounds on $s_9^{\standardtetrahedron}$. 

\section{Discussion}

Concerning the semidefinite relaxation described in Section~\ref{sec: SDP} there are still two open questions. First, for all investigated instances the optimal values of the second semidefinite program and third semidefinite program coincide. So the natural question that arises is
\begin{question}
    Are the semidefinite relaxations of Problem~\ref{prb: inner4} and Problem~\ref{prb: inner5} equivalent?
\end{question}
 Second, for all investigated instances the optimal values of the three semidefinite programs are integer which raises the question 
 \begin{question}
    Can the semidefinite relaxations of Problem~\ref{prb: inner3}, Problem~\ref{prb: inner4} and Problem~\ref{prb: inner5} be reduced to integer programs?
 \end{question}

In view of~\cite{GU}, a natural question is
\begin{question}\label{q:equivalence}
	How many non-equivalent maximal symplectic packings of the ball by equal balls exist ?
\end{question}
Note that, at the time of writing, the only answers to question~\ref{q:equivalence} are in dimension 4~\cite{McDuff1991}.
The ``natural'' definition of equivalence is: Two packings $\varphi,\psi:\amalg_kB(r)\hookrightarrow B(1)$ are equivalent if there exist a symplectomorphism $\phi:B(1)\to B(1)$ such that $\phi\circ\psi(\amalg_kB(r))=\varphi(\amalg_kB(r))$. For this notion, the answer to question~\ref{q:equivalence} is 1. Indeed,
\begin{Thm}[\cite{McDuff1991}]
	All isometric embeddings $\varphi,\psi:\amalg_kB^4(r)\hookrightarrow B^4(1)$ are isotopic through symplectic embeddings.
\end{Thm}

One can think of other notions of equivalence (for instance affine transformations mapping a packing onto another one or isometries,...) but since none seems more natural than the others we won't do it here.

\section{Acknowledgements}

We gratefully acknowledge the advice of Hansj\"org Geiges and Felix Schlenk. 
This research is part of a project in the SFB/TRR 191 `Symplectic Structures in Geometry, Algebra and Dynamics', funded by the DFG.

\bibliographystyle{alpha}
\bibliography{sympack_arxiv}

\medskip\noindent
Greta Fischer, Department of Mathematics and Computer Science, University of Cologne, Germany, 
\\\texttt{fischer@informatik.uni-koeln.de}

\medskip\noindent
Jean Gutt, Institute of Mathematics, University of Toulouse, France, 
\\\texttt{jean.gutt@math.univ-toulouse.fr}

\medskip\noindent
Michael J{\"u}nger, Department of Mathematics and Computer Science, University of Cologne, Germany, 
\\\texttt{juenger-sfb@informatik.uni-koeln.de}

\end{document}